\documentclass[10pt,a4paper,twoside,notitlepage]{article}

\usepackage{amsfonts}
\usepackage{amsmath}
\usepackage{amssymb}

\usepackage{graphicx}
\usepackage{psfrag}
\usepackage{color}
\usepackage{caption}
\usepackage{mathrsfs}
\usepackage[sans]{dsfont}

\newtheorem{thm}{Theorem}
\newtheorem{lem}{Lemma}[section]
\newtheorem{cor}{Corollary}[section]
\newtheorem{prop}{Proposition}[section]
\newtheorem{defn}{Definition}[section]
\newtheorem{rem}{Remark}[section]

\newcommand{\R}{\mathbb{R}}
\newcommand{\N}{\mathbb{N}}

\newcommand{\ve}{\varepsilon}
\newcommand{\n}{\noindent}

\newcommand{\vt}{\vartheta}
\newcommand{\s}{\sigma}

\newcommand{\mE}{\mathcal{E}}

\newcommand{\mO}{\mathcal{O}}

\begin{document}
   \begin{center}
      {\bf Multiplicative controllability for semilinear 
       reaction-diffusion equations\\ with 
      finitely many changes of sign
      \footnote{
This work was supported by the ÒIstituto Nazionale di Alta MatematicaÓ (INdAM), through the GNAMPA Research Project 2015: 
\lq\lq Analisi e controllo di equazioni a derivate parziali nonlineari''
(coordinator G. Floridia).
\\
Moreover, this research was 
performed in the framework of the GDRE CONEDP (European Research Group on \lq\lq Control of Partial Differential Equations'') issued by CNRS, INdAM and Universit\'e de Provence.
}}\\
      P. Cannarsa \\
      Department of Mathematics \\
      University of Rome ``Tor Vergata''\\
      00133 Rome, Italy \\
      email: cannarsa@mat.uniroma2.it \\
      and\\
      G. Floridia 
\footnote{The research of this author was carried out in the frame of Programme STAR, financially supported by UniNA and Compagnia di San Paolo.}
 \\
      Department of Mathematics and Applications \lq\lq R. Caccioppoli'' \\
      University of Naples \lq\lq Federico II''\\
Via Cintia, Monte S. Angelo
I-80126 Napoli, Italy\\
      email: 
      floridia.giuseppe@icloud.com\\
and \\
A.Y. Khapalov \\
Department of Mathematics\\
Washington State University,  Pullman, WA 99164-3113 USA; \\
fax: (1 509) 335 1188; tel. (1 509) 335 3172; 
e-mail: khapala@math.wsu.edu  
 \end{center}

\begin{abstract}
We study the global approximate controllability properties of a one dimensional semilinear reaction-diffusion equation governed via the coefficient of the reaction term.  It is assumed that both the initial and target states admit no more than  finitely many changes of sign. Our goal is to show that any target
state $ u^*\in H_0^1 (0,1)$,
 with as many changes of sign in the same order as the given initial data  $ u_0\in H^1_0(0,1)$,
 can be approximately reached in the $ L^2 (0,1)$-norm 
 at some time $T>0$. 
Our method employs shifting the points of sign change by making use of a finite sequence of 
initial-value pure diffusion problems.
\end{abstract}

\section{Introduction 
} 

Our main goal in this paper is to study the global approximate controllability properties of the following semilinear Dirichlet boundary value problem
\begin{equation}\label{1.1}
   \begin{cases}
\quad   u_t \; = \; u_{xx} \; + \; v(x,t)  u \; + \; f(u)
&\quad  {\rm in} \;\;\;Q_T = (0,1) \times (0, T)\,, \;\; 
\;\; T>0,
\\
\quad u (0,t) = u (1,t) = 0,
&\quad\qquad\quad t \in (0, T),
\\
\quad u\:\mid_{t = 0} \; = u_0\in H^1_0 ( 0,1 ). 
\end{cases}
\end{equation}
Here  $v \in L^\infty (Q_T) $ is a control function, which affects the reaction rate of the process described by \eqref{1.1}.
The nonlinear term $f:\R\rightarrow\R$ is assumed to be a Lipschitz function satisfying $f(0)=0$.

Let us recall that, in  general terms, an evolution system  is called globally approximately
controllable in a given space $H$ at time
$ T > 0$, if any initial state in $ H$ can be steered  
into any neighborhood of any desirable target state  at time $ T$,
by selecting a suitable control. 

Historically, the concept of controllability  emerged  in the context of linear ordinary differenetial equations and was motivated by numerous engineering applications. 
Then it  was extended to various linear partial differential equations governed by {\it additive} locally distributed (i.e., supported on a bounded subdomain of the space 
domain), lumped (acting at a point), and boundary controls (see Fattorini in \cite{Fatt}, Fattorini and Russell in \cite{FattRuss}, and Zabczyk in \cite{Zab}). 
Methodologically, these studies were typically based on the {\em linear} duality pairing technique between  the  control-to-state mapping at hand and its dual observation map.
When  this mapping is nonlinear, as it happens in \eqref{1.1},  the aforementioned  approach does not apply  and the above-stated concept of controllability becomes, in general, unachievable. 
For example, the above control system cannot be steered to any nonzero target starting from the origin $u_0=0$ (see more about that in  Remark \ref{r1.2} below).

It is well-known that \eqref{1.1} 
can be linked to  various applied reaction-diffusion models such as chemical reactions, nuclear chain reactions, and biomedical models (see \cite{KhBook}, and the references therein). 
More generally, reaction-diffusion equations or systems 
describe how the concentration of one or more substances 
changes under the influence of some processes such as local chemical reactions, where
substances are transformed into each other, and diffusion which causes substances to spread out 
in space.
 
 Unfortunately, additive controls (see, e.g., \cite{ACF} and \cite{CMV1}) are unfit to treat such problems because, for example, they would require inputs with high energy levels or they are not available
due to the physical nature of the process at hand. On the other hand, an approach based on multiplicative controls, where the coefficient $v$ in \eqref{1.1} is used to change the main physical characteristics of the system  at hand, seems realistic.
 
In the area of multiplicative controllability for partial differential equations we would like to mention the pioneering work \cite{BMS} by Ball, Marsden and Slemrod  establishing the approximate controllability of the rod  and  wave equations, based on the implicit nonharmonic Fourier series approach. 
 Further principal contributions to the field of  the multiplicative controllability of linear and semilinear parabolic and hyperbolic equations and of a number of swimming models were made  by Khapalov, see  \cite{Kh1}-\cite{KhBook} and the references therein.
Substantial progress has also been made in the study of the controllability properties of
the Schr\"odinger equation by Beauchard, Coron, Boscain et al., Ervedoza and Puel, Nersesyan, see 
 \cite{B1}, \cite{BC}, \cite{CMSB}, \cite{EP} and \cite{Ner} and  the references therein. Let us also mention along these lines  the work by  Beauchard \cite{B2} for the beam model,  the works by Lin et al. \cite{Lin1} and  \cite{Lin2}, and the work by  Fernandez and Khapalov \cite{FKh} on the bilinear controllability of parabolic equations.

 In \cite{Kh1}, Khapalov studied the global nonnegative approximate controllability of the one dimensional 
 nondegenerate semilinear convection-diffusion-reaction equation governed in a bounded domain via bilinear control. Similar results were obtained for degenerate parabolic equations by Cannarsa and Floridia in \cite{CF1},\cite{CF2},\cite{F1} and \cite{F2}.
In  \cite{CanKh} Cannarsa  and Khapalov established  an approximate controllability property  for nondegenerate linear equations in suitable classes of functions that change sign.

%
 In this paper, we are interested in the multiplicative controllability of the semilinear reaction-diffusion system \eqref{1.1} when both the initial and target states admit a finite number of points of sign change. This fact introduces substantial differences with respect to the above works. Indeed, on the one hand, we manage to extend all the approximate controllability results of \cite{CanKh} to semilinear equations as well as those of \cite{Kh1} to initial/target states that may change sign. On the other hand, here we introduce a new technique of proof. 
 
 In \cite{CanKh}, an implicit \lq\lq {\it continuation argument}'' was employed to justify the fact that one can always continue to move the points of sign change until their target positions have been reached. Such a technique was mainly qualitative but sufficient to obtain the conclusion due to the linear structure of the equation at hand. In this paper, on the contrary, a more quantitative approach is needed because the equation of interest is nonlinear. 
 
 Indeed,  we  give an explicit construction of the controls required for the steering process. Such controls are, essentially, obtained by splitting $[0,T]$ into finitely many time intervals 
$$[0,T]=[0,S_1]\cup[S_1,T_1]\cup\cdots\cup[T_{N-1},S_N]\cup[S_{N},T_N]\cup[T_N,T]$$
on which two alternative actions are applied: on $[S_k,T_k]$ we choose suitable initial data, $w_k$, in pure diffusion problems ($v\equiv 0$) to move the points of sign change to their desired location, whereas on $[T_{k-1},S_k]$ we use piecewise static multiplicative controls $v_k$ to attain such $w_k$'s as intermediate final conditions. More precisely, on $[S_k,T_k]$ we make use of the  boundary problems 
\begin{equation*}
   \begin{cases}
\quad   w_t \; = \; w_{xx} \;  + \; f(w), &
 {\rm in} \;\;\;(0,1) \times [S_{k}, T_{k}],
\\
\quad w (0,t) = w (1,t) = 0,
& t \in [S_k, T_{k}],
\\
\quad w\mid_{t = S_k} \; = w_{k}(x),
& w_{k}'' (x) \mid_{x = 0, 1} = 0,
\end{cases}
\end{equation*}
where the $ w_{k}$'s are viewed as  control parameters to be chosen to generate suitable curves of sign change, which have to be continued along all the $N$ time intervals $[S_k,T_k]$ until each point has reached the desired final position. 
In order to fill the gaps between two successive $[S_k, T_{k}]$'s, on $[T_{k-1},S_k]$ we construct $v_k$ that steers the solution of 
\begin{equation*}
\begin{cases}
\quad   u_t \; = \; u_{xx} \; + \; v_k(x,t)  u \; + \; f(u),
&\quad  {\rm in} \;\;\;
 (0, 1) \times [T_{k-1}, S_{k}],  
\\
\quad u (0,t) = u (1,t) = 0,
&\quad\quad\qquad\;\; t \in [T_{k-1}, S_{k}],
\\
\quad u\:\mid_{t = T_{k-1}} \; = u_{k-1}+r_{k-1} ,
\end{cases}
\end{equation*}
 from $u_{k-1}+r_{k-1}$ 
to  $w_k,$
where $ u_{k-1}$ and $ w_{k}$ have  
the same points of sign change, 
 and  $\|r_{k-1}\|_{L^2 (0, 1)}$ is small.
The fact that such a process can be completed within a finite number of steps is an important point of the proof.
It follows from precise estimates that guarantee that the sum of the distances of each branch of the null set of the resulting solution of \eqref{1.1} from its target points of sign change decreases  at a linear-in-time rate for curves which are still far away from their corresponding target points, while the  error caused by the possible displacement of points already near their targets is negligible. 
 
 We believe that the techniques of this paper could be useful to study multidimensional diffusive systems on  Riemannian manifolds of low dimension or special structure.
 
 \subsubsection*{Motivations and future perspectives}
In this section we present some applications. 
The following 
nuclear model   
is a typical example of the applications we plan to study, that is, a reaction-diffusion model in a fissionable material (see Section 2.7 of \cite{Salsa}).
 By shooting neutrons into a uranium nucleus it may happen that the nucleus
breaks into two parts, releasing other neutrons already present in the nucleus and
causing a chain reaction. At a macroscopic level, the free neutrons diffuse like a chemical in a porous medium, where reaction and diffusion are competing.
 Some macroscopic aspects of this phenomenon can be described by means of a simplified reaction-diffusion model, like in \eqref{1.1}, where $u$ is the neutron density and the multiplicative coefficient $v$ is the fission rate.
 
 Our study of the reaction-diffusion models is also motivated 
 by mathematical models of tumor growth  (see, e.g., Friedman in \cite{Fr bio} and Perthame in \cite{Perthame}). 
There are three distinct main stages in the growth of a tumor (see \cite{Roose}, \cite{Tosin} and \cite{Anderson}) before it becomes so large that it causes patients to die or reduces permanently their quality of life: avascular (tumors without blood vessels), vascular, and metastatic. 
From a clinical point of view, vascular and metastatic tumor growth
are what cause the patient to die. So modeling and understanding
these processes is crucial for cancer therapy. Nevertheless, avascular tumor growth is much simpler to model mathematically, and yet contains many of the phenomena 
that one needs to address also in a general model of vascular or metastatic tumor growth.
In the review \cite{Roose}, Roose, Chapman and Maini describe a continuum mathematical model of avascular tumor growth.  This model consist of reaction-diffusion-convection equations and was introduced by Casciari, Sotirchos and Sutherland in \cite{Casciari}.
Mathematical models describing continuum tumor cell populations and their development classically consider the interactions between the cell number density and one or more chemical species that provide nutrients or influence the cell cycle events. 
The model introduced in \cite{Casciari} 
can reduce to the following simplified 
system
$$\frac{\partial u_{i}}{\partial t} \; = \frac{\partial^2 u_{i}}{\partial^2x} \; + \; v_i(x,t)  u_i \; + \; f(u_i),\;\;\; i=1,\ldots,N\;\; (N\in\N),$$
where $u_i$ are the concentrations of the chemical species and $v_i$ are the net rate of consumption/production of the chemical
species both by the tumor cells and due to chemical reactions with other species. 
 In future works we will address controllability issues for systems of
reaction-diffusion equations of the type outlined 
in the survey paper \cite{BP00} 
(see Section 7, {\it \lq\lq Control problems''}).
\subsubsection*{Outline of the paper}
In Section 2, we give the precise formulation of the problem and state our approximate controllability result for system \eqref{1.1} (Theorem~\ref{th1.1}) together with some of its consequences. \\
In Section 3, we explain the structure of the proof and, admitting  two essential technical tools (Theorem~\ref{th2.2} and Theorem~\ref{th preserving}), we proceed with the proof of Theorem~\ref{th1.1}.\\
Section 4 deals with the proof of Theorem~\ref{th2.2}: a controllability result for pure diffusion problems. \\
Finally, 
Section 5 is devoted to the proof of Theorem~\ref{th preserving}: a smoothing result intended to attain suitable intermediate data while preserving the already-reached points of sign change.

\section{Problem formulation and results} 
Our main goal in this paper is to study the global approximate controllability properties of the semilinear Dirichlet boundary value problem \eqref{1.1}
\begin{equation*}
   \begin{cases}
\quad   u_t \; = \; u_{xx} \; + \; v(x,t)  u \; + \; f(u)
&\quad  {\rm in} \;\;\;\,Q_T = (0,1) \times (0, T)\,, \;\; 
\;\; T>0,
\\
\quad u (0,t) = u (1,t) = 0,
&\quad\qquad\quad t \in (0, T),
\\
\quad u\:\mid_{t = 0} \; = u_0\in H^1_0 ( 0,1 ). 
\end{cases}
\end{equation*}
Here  $v \in L^\infty (Q_T) $ is a bilinear control. 
The nonlinear term $f:\R\rightarrow\R$ is supposed to be a Lipschitz function with $f(0)=0,$ differentiable at $0$ and $L$ will denote a Lipschitz constant for $f,$ that is,
\begin{equation}\label{1.2a}
\mid f(u') - f (u) \mid \; \leq \; L \mid u' - u \mid ,\;\;\;\; \forall u,u' \in \R.
\end{equation}
\begin{rem}
\label{r1.2}\rm
We note that system \eqref{1.1} cannot be steered anywhere from the origin.
Moreover, if $u_0 (x) \geq 0$ in $(0,1)$,  then the strong maximum principle 
demands that the respective solution to \eqref{1.1} remains nonnegative  at any moment of time, regardless of the choice of $ v$. This means that system \eqref{1.1} cannot be steered from any such $ u_0 $ to any target state which is negative on a nonzero measure set in the space domain. We remark that the strong maximum principle for linear parabolic PDEs (see, e.g, Chapter 2 in \cite{Fr}, p. 34) can be extend to semilinear parabolic system \eqref{1.1}. Indeed, since $f(0)=0$ and $f(u)$ is differentiable at $0,$ 
the term $ f(u (x,t))$ in $\eqref{1.1}$ can be represented as $ \beta (x,t) u (x,t) $,  where $ \beta = \dfrac{f(u)}{u} \in L^\infty (Q_T).$
\end{rem}

Let us start 
with the well-posedness for the system \eqref{1.1}.
\subsubsection*{Functional setting and Well-posedness}
Hereafter, we use the standard notation for Sobolev spaces, in particular,   
\begin{align*}
H^1 (0,1) &= \{\phi\in L^2 (0,1) \mid \; \phi_{x} \in L^2 (0,1)
\}\\
H_0^1 (0,1) &= \{\phi\in H^1 (0,1)  \mid 
 \phi(0)=\phi(1)=0 
\}\\
H^2 (0,1) &= \{\phi\in H^1 (0,1)  \mid \; 
\phi_{xx}\in L^2 (0,1) 
\}.
\end{align*}
By classical well-posedness results (see, for instance, Theorem 6.1 in \cite{Lad}, pp. 466-467) problem (\ref{1.1}) 
with initial data $u_0\in L^2(-1,1)$ admits a unique solution 
 $$u\in L^2(0,T; H^1_0 (0,1) )\cap  C ( [0, T]; L^2 (0,1)).$$ 
Furthermore, if $u_0\in H^1_0(0,1)$, then the solution $u$ of problem (\ref{1.1}) satisfies
$$u\in H^1(0,T; L^2 (0,1) )\cap  C ( [0, T];  H^1_0(0,1))\cap  L^2 ( 0, T;  H^2(0,1)).$$
\subsubsection*{Problem formulation}
In this paper,  we assume that $u_0\in H_0^1 (0,1)$ has  {\it finitely many zeros},
 that is, there exist points 
 $$0=x_0^0<x^0_1<\cdots<x^0_n<x_{n+1}^0=1$$ such that
  $$u_0(x)=0\Longleftrightarrow x=x_l^0,\; l=0,\ldots,n+1.$$ 
  Moreover, we assume that the interior zeros ($x_l^0,\, l=1,\ldots,n$) are {\it points of sign change}, that is, for $l=1,\ldots,n,$ 
  $$u_0(x)u_0(y)<0,\;\;\; \forall x\in \left(x^0_{l-1},x^0_l\right),\, \forall y\in \left(x^0_{l},x^0_{l+1}\right).
   $$
We will refer to such functions $u_0$ as the ones with {\it finitely many changes of sign}. \\

Our goal is to show that any target $ u^* \in H_0^1 (0,1)$, with as many changes of sign {\it in the same order} as the given $ u_0$, can be approximately  reached in the $ L^2 (0,1)$-norm   at some time $ T>0$. 
By the above expression we mean that, denoting by $x_l^*,\, l=0,\ldots,n+1$, the zeros
 of $u^*,$ we have
 $$u_0(x)u^*(y)>0,\;\;\; \forall x\in \left(x^0_{l-1},x^0_l\right),\, \forall y\in \left(x^*_{l-1},x_{l}^*\right), \text{ for } l=1,\ldots,n+1.$$

\begin{rem}\label{optimal}\rm
The matching of the initial and target states  
seems optimal. Indeed, the strong maximum principle (see Remark \ref{r1.2}), applied on the subdomain of $Q_T$  
delimited by any two adjacent curves of sign change (in the sense of Definition \ref{scgs} and Lemma \ref{p1}.),  
prevents the appearance of new zeros of $u(x,t)$ within this area. 
\end{rem}
Let us start with the following definition.
\begin{defn}\label{static}
We say that a 
function $v\in L^\infty(Q_T)$ is {\it piecewise static}, if there exist $m\in\N,$ $c_k(x)\in L^\infty(0,1)$ and $t_k\in [0,T], \,t_{k-1}<t_k,\, k=1,\dots,m$ with $t_0=0 \mbox{ and } t_m=T,$ 
such that $$v(x,t)=c_1(x)\mathds{1}_{[t_{0},t_1]}(t)+\sum_{k=2}^m c_k(x)\mathds{1}_{(t_{k-1},t_k]}(t),$$ where $\mathds{1}_{[t_{0},t_1]}\,  \mbox{  and  }  \,\mathds{1}_{(t_{k-1},t_k]}$ are the indicator function of $[t_{0},t_1]$ and $(t_{k-1},t_k]$, respectively.
\end{defn}
Here we present our main result for system (\ref{1.1}).
\begin{thm}\label{th1.1} Let $ u_0 \in  H_0^1 (0,1)$. Assume that $ u_0$ has finitely many points of sign change.
 Consider any  $ u^*\in H_0^1 (0,1)$ which has exactly as many points of sign change in the same order as $ u_0 $. Then, for any $ \eta > 0$ there are a $ T=T(\eta,u_0,u^*) > 0$ and  a piecewise static multiplicative control $ v=v(\eta,u_0,u^*) \in L^\infty (Q_T) $ such that for the respective solution $u$ to $\eqref{1.1}$ the following inequality holds
\begin{equation*}
\| u(\cdot, T) - u^* \|_{L^2 (0,1)} \; \leq \; \eta.
\end{equation*}
\end{thm}

\vspace{-1.0cm}
\setlength{\unitlength}{1mm}
\begin{picture}(150,55)(0,0)
\linethickness{1pt}
\put (2,15){\vector(1,0){85}}
\put (10,10){\vector(0,1){35}}
\qbezier(10,15)(16,32)(23,18)
\qbezier(23,18)(34,1)(45,10)
\qbezier(45,10)(65,25)(70,15)
{
\qbezier(10,15)(19,42)(27,19)
\qbezier(27,19)(34,-15)(48,18)
\qbezier(48,18)(65,55)(70,15)
}
\put(23,10){\makebox(0,0)[b]{$x^0_{1}$}}
\put(54,10){\makebox(0,0)[b]{$x^0_{2}$}}
\put(31,15){\makebox(0,0)[b]{$x^*_{1}$}}
\put(44,15){\makebox(0,0)[b]{$x^*_{2}$}}
\put(64,20){\makebox(0,0)[b]{$u_0$}}
\put(62,37){\makebox(0,0)[b]{$u^*$}}
\put (25,15){\vector(1,0){5}}
\put (25,15){\vector(0,1){5}}
\put (25,14){\makebox(0,0)[b]{$\bullet$}}
\put (52.7,15){\vector(-1,0){5}}
\put (52.7,15){\vector(0,1){5}}
\put (52.7,14){\makebox(0,0)[b]{$\bullet$}}

\put(8,40){\makebox(0,0)[b]{$u$}}
\put(85,10){\makebox(0,0)[b]{$x$}}

\put(8,10){\makebox(0,0)[b]{$0$}}
\put(70,15){\line(0,1){2}}
\put(70,10){\makebox(0,0)[b]{$1$}}
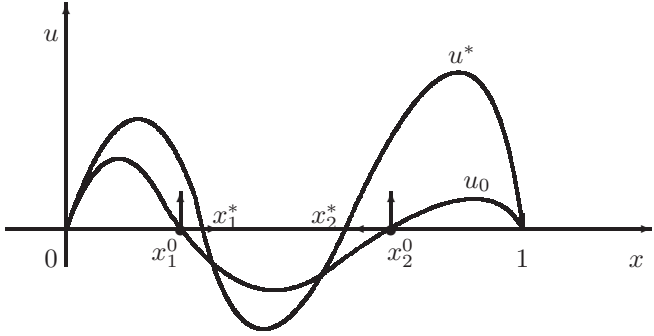
\put(35,-5){\makebox(0,10)[b]{ Figure 1. \; \rm Control of 
two points of sign change.}}
\end{picture}


\subsection{Further results}
We present in this section two results that generalize Theorem \ref{1.1} and are easy consequences of such theorem.
\begin{cor} \label{th1.2} 
Let $ u_0,\,u^* \in  H_0^1 (0,1).$ 
Assume that $ u_0$ and $u^*$ have finitely many points of sign change and
 the amount of points of sign change of $u^*$ is less than the one
  of $u_0.$
 Then, for any $ \eta > 0$ there are a $ T=T(\eta,u_0,u^*) > 0$ and  a piecewise static multiplicative control $ v=v(\eta,u_0,u^*) \in L^\infty (Q_T) $ such that for the solution $u$ to $\eqref{1.1}$ the following inequality holds
\begin{equation*}
\parallel u(\cdot, T) - u^* \parallel_{L^2 (0,1)} \; \leq \; \eta.
\end{equation*}
\end{cor}
In the following Remark \ref{rem th1.2} we clarify the statement of Corollary \ref{th1.2}.
\begin{rem}\label{rem th1.2}\rm
We explain the statement of Corollary \ref{th1.2} by the following example. Let us consider an interval $ (0, x^0_1)$ of positive values of $u_0$ followed by an interval $ (x^0_1, x^0_2)$ of negative values of $ u_0 (x)$, which in turn is followed by  an interval $ (x^0_2, x^0_3)$ of positive values of $ u_0 (x)$ and so forth. Then the merging of the respective two points of sign change $ x^0_1 $ and  $ x^0_2$ will result in one single  interval $ (0, x^0_3)$ of positive 
values. 
\end{rem}
In the following picture we describe the situation discussed in
 Remark \ref{rem th1.2} in the particular case $0=x_0^0<x_1^0<x_2^0<x_3^0=1.$

\vspace{-0.7cm}
\setlength{\unitlength}{1mm}
\begin{picture}(150,55)(0,0)
\linethickness{1pt}
\put (2,15){\vector(1,0){85}}
\put (10,10){\vector(0,1){35}}
\qbezier(10,15)(16,32)(23,18)
\qbezier(23,18)(34,1)(45,10)
\qbezier(45,10)(65,25)(70,15)
{
\qbezier(10,15)(65,25)(70,15)}
\put(25,10){\makebox(0,0)[b]{$x^0_{1}$}}
\put(56,10){\makebox(0,0)[b]{$x^0_{2}$}}
\put(36,8){\makebox(0,0)[b]{$u_{0}$}}
\put(40,20){\makebox(0,0)[b]{$u^*$}}

\put(8,40){\makebox(0,0)[b]{$u$}}
\put(85,10){\makebox(0,0)[b]{$x$}}

\put(8,10){\makebox(0,0)[b]{$0$}}
\put(70,15){\line(0,1){2}}
\put(70,10){\makebox(0,0)[b]{$1$}}
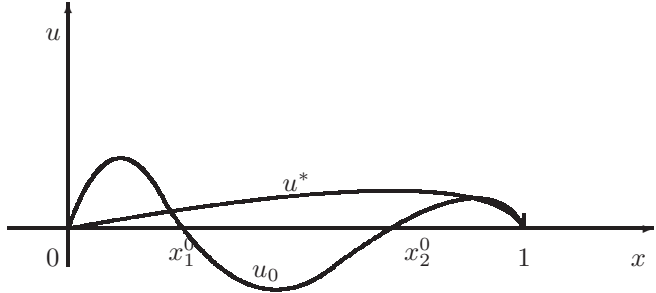
\put(42,0){\makebox(0,10)[b]{{ Figure 2.} \;$u_0, {
u^*}$:\;{\rm 
merging of the points of change of sign.} }}
\end{picture}
{\bf Proof (of Corollary \ref{th1.2}).}
Corollary \ref{th1.2} follows from Theorem \ref{th1.1}. Indeed, all the target states described in Corollary \ref{th1.2} can be approximated in $L^2(0,1)$ by those in Theorem \ref{th1.1}.\\

It may be worth noting that the following generalized approximate controllability property can be deduced from  Corollary \ref{th1.2}.
\begin{cor}\label{cor CK}
 Let $u_0$ and $u^*$ be given  in $ L^2 (0,1)$.  Then, for any $ \eta > 0$ there exists $u^\eta_0 \in H_0^1 (0,1)$ such that $\parallel u^\eta_0 - u_0 \parallel_{L^2 (0,1)} \; < \; \eta,$ and there exist $ T=T(\eta,u_0,u^*)> 0$ and a piecewise static multiplicative control $v=v(\eta,u_0,u^*)\in L^\infty(Q_T)$ such that the solution $u$ to
$$
   \begin{cases}
\quad  u_t \; = \; u_{xx} \; + \; v(x,t)  u \; + \; f(u)
&\quad {\rm in} \;\;\;
   Q_T=(0,1)\times (0,T)
\\
\quad u (0,t) = u (1,t) = 0
&\quad \quad\qquad t \in (0, T)
\\
\quad u(\cdot, 0)=u^\eta_0 \in H_0^1 (0,1)
\end{cases}
$$ 
satisfy the following inequality
$$
\| u(\cdot, T) - u^* \|_{L^2 (0,1)} \; \leq\eta.
$$
\end{cor}

The proof of Corollary \ref{cor CK} is 
similar to one of the Corollary 2 of \cite{CanKh}, which deals with 
linear systems.

\section{
Control Strategy for the proof of the main result}\label{Sec Control Strategy} 
In this section, we fix the notation and we
 introduce a control strategy  
to obtain the complete proof of Theorem \ref{th1.1} 
 in Section \ref{proof main result}.\\
 
Let us fix a number $\vt\in (0,1)$ to
 be used in whole the paper.\\

\noindent Now, we recall some useful functional spaces. 
\subsubsection*{H\"older continuous spaces}
We define the H\"older spaces
$$C^{\vt} ([0,1]):=\left\{w\in C ([0,1])\,:\:\sup_{x,y\in[0,1]}\frac{|w(x)-w(y)|}{|x-y|^\vt}
<+\infty\right\},$$
$$C^{2+\vt} ([0,1]):=\left\{w\in C^2 ([0,1])\,:\:w''\in C^{\vt} ([0,1])
\right\}.$$
Let $Q_T=(0,1)\times(0,T)$. Let us define the following 
 spaces of time dependent functions 
  $$C^{\vt, \frac{\vt}{2}} (\overline{Q}_{T}):=\left\{u\in C (\overline{Q}_{T})\,:\:\sup_{x,y\in[0,1]}\frac{|u(x,t)-u(y,t)|}{|x-y|^\vt}+\sup_{t,s\in [0,T]}\frac{|u(x,t)-u(x,s)|}{|t-s|^\frac{\vt}{2}}<+\infty\right\},
 $$
 $$C^{2, 1} (\overline{Q}_{T}):=\left\{u:\overline{Q}_{T}\longrightarrow \R\,:\exists\:u_{xx},u_x,u_t\in C(\overline{Q}_{T})
    \right\},
 $$  
  $$C^{2+\vt, 1+\frac{\vt}{2}} (\overline{Q}_{T}):=\left\{u\in C^{2,1} (\overline{Q}_{T})\,:\:u_{xx},u_x,u_t\in C^{\vt, \frac{\vt}{2}} (\overline{Q}_{T})
    \right\}.
 $$ 
\subsubsection*{Notation}
Given $N\in\N,$ let us set 
$\R^N_+=\{(a_1,\ldots,a_N)\,|\,a_k\in\R,\,a_k>0, \,k=1,\ldots,N\}$.
For every $$
(\tau_1,\ldots,\tau_N)=(\tau_k)_1^N\in\R^N_+,\;\;
(\s_1,\ldots,\s_N)=(\s_k)_1^N\in\R^N_+,
$$
we define 
\begin{equation}\label{N1}
T_0:=0,\qquad S_k:=T_{k-1}+\s_k, \qquad T_k:=S_k+\tau_k,\qquad\quad k=1,\ldots,N.  \qquad(\footnote{We note that $T_0=0,\, S_1=\s_1,\, T_1=\s_1+\tau_1,$ and $\displaystyle S_k=\sum_{h=1}^{k-1}(\s_{h}+\tau_{h})+\s_k,\; T_k=\sum_{h=1}^{k}(\s_h+\tau_h),$ $\forall k=2,\ldots,N.$})
\end{equation} 
Noting that $0=T_0<S_{k}<T_{k}\leq T_N,\;k=1,\ldots,N,$ we consider 
the following partition of $[0,T_N]$ in $2 N$ intervals:
\begin{equation}\label{N2}
[0,T_N]=[0,S_1]\cup[S_1,T_1]\cup\cdots\cup[T_{N-1},S_N]\cup[S_{N},T_N]=\bigcup_{k=1}^N\left(\mO_k \cup \mE_k\right),
\end{equation}
where, for every $k=1,\ldots,N,$ 
we have set $\mO_k:=[T_{k-1},S_k]$ 
and $\mE_k:=[S_k,T_k].$ 
\subsubsection*{Outline and main ideas for the proof of Theorem \ref{th1.1}}
We obtain the proof of Theorem \ref{th1.1} in Section \ref{proof main result} using the partition introduced in \eqref{N2} and applying two alternative control actions: 
on $[S_k,T_k]$ we choose suitable initial data, $w_k$, in pure diffusion problems ($v\equiv 0$) to move the points of sign change to their desired location (Section \ref{pure diffusion}), whereas on $[T_{k-1},S_k]$ we give a smoothing result to preserve the reached points of sign change and attain such $w_k$'s as intermediate final conditions,
using piecewise static multiplicative controls $v_k$ (Section \ref{sec smoth}). 
 \subsection{Controllability for initial-value 
pure diffusion problems on disjoint time intervals
 }\label{pure diffusion} 
In this section, we outline the main idea of the proof of Theorem \ref{1.1}.\\

Let $N\in\N.$ For any fixed $
(\s_1,\ldots,\s_N) 
\in \R^{N}_+,$ 
let us consider a generic $(\tau_1,\ldots,\tau_N)
\in \R^{N}_+$ and, for
$ k = 1, \ldots, N,$ recalling 
\eqref{N1}-\eqref{N2}, 
 let us introduce the following initial value 
pure diffusion problems on disjoint time intervals 
\begin{equation}\label{2.1}
   \begin{cases}
\quad   w_t \; = \; w_{xx} \;  + \; f(w), \;\;\;\;
 {\rm in} \;\;\;Q_{\mE_k
 } = (0,1) \times [S_{k}, T_{k}
 ],
\\
\quad w (0,t) = w (1,t) = 0,
\qquad\qquad\qquad\;\;\:\;\;\, t \in [S_k, T_{k}],
\\
\quad w\mid_{t = S_k} \; = w_{k}(x),
\;\;\;\;\qquad w_{k}'' (x) \mid_{x = 0, 1} = 0. 
\end{cases}
 \;(\footnote{$w_{k}'' (x) \mid_{x = 0, 1} = 0$ is a compatibility condition, see \cite{Lad}, pp. 452-453, where it is introduced in a more general parabolic problem.})
\end{equation}
We will consider the initial data $w_{k} $ and times 
 $\tau_k
 $,\, $ k = 1, \ldots, N,$ 
as control parameters, where $ w_{k} $'s belong to $ C^{2+ \vt} ([0,1]),$ 
with $\vt \in(0, 1)$ fixed at the beginning of this section.
 \begin{rem}\label{solution of dtp}
 \rm
For every $k=1,\ldots,N,$ by a  classical well-posedness result (see \cite{Lad}, pp. 452-453), if $ w_{k}\in C^{2+ \vt} ([0,1]),$ then any initial-value problem in \eqref{2.1} has a unique classical solution $W_k(x,t)$ on $\overline{Q}_{\mE_k}$
 and
 $ W_k \in C^{2+\vt, 1+\vt/2} (\overline{Q}_{\mE_k}
 ).$
 \end{rem}
\begin{defn}\label{def ds}
We call 
solution of \eqref{2.1} the function defined in $\displaystyle (0,1)\times {\small{\small \bigcup_{k=1}^N}}[S_k,T_k]
$ as
$$w(x,t)=W_k(x,t),\;\;\; \forall (x,t)\in(0,1)\times [S_{k}, T_{k}
],\;\;k=1,\ldots,N,$$
where $W_k,$ for every $k=1,\ldots,N,$ is the unique solution on $(0,1)\times [S_{k}, T_{k}
]$ of the $k^{th}$ problem in \eqref{2.1}, with initial state $w_k.$
\end{defn}
\begin{rem}\label{ds}\rm
We observe that a 
 solution of \eqref{2.1} is a collection of solutions of finitely many problems which are set on disjoint time intervals. Therefore, it is 
 independent of the choice of $(\s_k)^N_1.$ We prefer to give Definition \ref{defTID} for a fixed $(\s_k)^N_1,$ just for technical purposes that will be clarified in the sequel (see Theorem \ref{th2.2} and Section \ref{proof main result}).
\end{rem}
\begin{defn}\label{defTID}
Let $u_0\in H^1_0(0,1)$ be a function with 
 $n$ points of sign change. For every fixed $N\in\N$ and $(\s_k)^N_1\in\R^N_+,$ we call a finite \lq\lq family of Times and Initial Data'' 
of 
\eqref{2.1} associated to $u_0,$ a set of the form
  $\displaystyle\left\{(\tau_k)_1^{N},(w_{k})_1^{N}\right\}
  $ 
  such that 
\begin{description}
\item[$\star$] \quad$(\tau_k)_1^N\in\R^N_+;$
\item[$\star$] $\;\;\;\;w_{k}\in C^{2+ \vt} ([0,1]),$ for all $k = 1, \ldots, N,$ satisfies:
\begin{enumerate}
\item 
$w_k(0)=w_k(1)=0,\;w''_k(0)=w''_k(1)=0;$ 
\item 
$w_1$ and $u_0$ have the same points of sign change, in the same order as the points of sign change of $u_0;$
\item for $k = 2, \ldots, N,$ $ w_{k} (\cdot) $ and $w(\cdot, T_{k-1}
)$ have the same points as the points of sign change, in the same order of sign change of $u_0,$ where $w$ is the solution of \eqref{2.1}.
\end{enumerate}
\end{description}
\end{defn}
\begin{thm}\label{th2.2} 
Let $ u_{0}\in 
H^1_0(0,1)$ have $n$ points of sign change at $ x^0_{l} \in (0,1),\, 
$ 
with $0=:x^0_0<x^0_l<x^0_{l+1}\leq x^0_{n+1}:=1,\;\; l=1,\ldots,n$.
Let $ x^*_{l} \in (0,1),\,  l=1,\ldots,n,$ be such that $0:=x^*_0<x^*_l<x^*_{l+1}\leq x^*_{n+1}:=1.$ 
Then, for every $\ve>0$ there exist $N_{\ve}\in\N$ and a finite family of times and initial data 
$ \left\{(\tau_k)^{N_{\ve}}_1, (w_{k})^{N_{\ve}}_1\right\}$ 
 such that, for any $(\s_k)^{N_\ve}_1\in\R^N_+,$   
the 
solution $w^{\ve}$ of 
problem \eqref{2.1} satisfies
$$ w^{\ve} (x, T_{N_{\ve}}
)=0 \qquad \Longleftrightarrow \qquad x=x_l^{\ve},\;\;l=0,\ldots,n+1,$$
for some points $x_l^{\ve}\in (0,1),\; 0:=x^{\ve}_0<x^{\ve}_l<x^{\ve}_{l+1}\leq x^{\ve}_{n+1}:=1,\;l=1,\ldots,n,$ such that
$$\sum_{l=1}^n|x^*_l-x^{\ve}_l|<\ve.$$
Moreover, $w^\ve(\cdot, T_{N_{\ve}}
)$ has the same order of sign change as $u_0.$ 
\end{thm}

This theorem is obtained in Section \ref{ADCP} by proving a series of preliminar results.

\subsection{Smoothing result to preserve the reached points of sign change and to obtain the intermediate data $w_k$}\label{sec smoth}
 In this section we introduce a smoothing result to preserve the reached points of sign change and attain regular intermediate final conditions $w_k$'s. \\
 
Let $N\in\N.$ For any fixed $
(\tau_1,\ldots,\tau_N) 
\in \R^{N}_+,$ 
let us consider a generic $(\s_1,\ldots,\s_N)
\in \R^{N}_+$ and, for
$ k = 1, \ldots, N,$ recalling 
\eqref{N1}-\eqref{N2}, 
given $u_{k-1},\,r_{k-1}\in H^1_0(0,1),\,v_k\in  L^\infty((0,1)\times  [T_{k-1},S_k]
),$ let us introduce the following problem
\begin{equation}\label{u pres}
   \begin{cases}
\quad   u_t \; = \; u_{xx} \; + \; v_k(x,t)  u \; + \; f(u)
&\quad  {\rm in} \;\;\;\,Q_{\mO_k}= (0, 1) \times [T_{k-1}, T_{k-1}+\s_k],  
\\
\quad u (0,t) = u (1,t) = 0,
&\qquad\qquad\qquad\qquad t \in [T_{k-1}, T_{k-1}+\s_k],
\\
\quad u\:\mid_{t = T_{k-1}} \; = u_{k-1} + r_{k-1}\in H^1_0(0,1).
\end{cases}
\quad(\footnote{We recall that $S_k=T_{k-1}+\s_k$.})\;
\end{equation}

Our goal is to show that, given $T_{k-1},$ there exists a sufficiently small $\s_k>0$ such that, provided  $ \|r_{k-1}\|_{L^2 (0, 1)} $ is small, we can steer the system (\ref{u pres}) from $u_{k-1} + r_{k-1}$ 
to a neighborhood of any state $w_k,$ 
where $  
w_k$  
has the same $n$ changes of sign as $u_{k-1}$,  
in the same order of sign change. In Section \ref{proof main result} we will apply the following result 
 to obtain regular $w_k$'s, that satisfy suitable properties.

\begin{thm}\label{th preserving} 
Let $ u_{k-1}, r_{k-1}, w_k\in 
H^1_0(0,1).$ Let $u_{k-1}$ and $w_{k}$ have the same $n$ points of sign change
in the same order of sign change.
Then, for every $\eta>0$ there exist a sufficiently small $\s_k=\s_k(\eta, u_{k-1}, 
 w_k)>0 
$
 and  
a piecewise static bilinear control $v_k=v_k(\eta, u_{k-1}, 
  w_k)\in L^\infty((0,1)\times(T_{k-1},S_k))$ 
such that
$$\|U_k (\cdot, S_k) - w_{k}(\cdot)\|_{L^2(0,1)}\leq\eta+ C_k
\|r_{k-1}\|_{L^2(0,1)}, 
$$ 
 where 
 $U_k$ is the solution of \eqref{u pres} on $(0,1)\times[T_{k-1},S_k]$ and 
 $C_k=C(u_{k-1}, w_k
 )$ is a positive constant. 
\end{thm}
 
 The above theorem 
 is proved in Section \ref{regular section}.
 \subsection{Proof of 
 Theorem \ref{th1.1}
 }\label{proof main result}
 Let us start this section by the following Lemma.
  \begin{lem}\label{pre proof th}
 Let $u_0\in H^1_0(0,1)$ be a function with 
 $n$ points of sign change. Let $\displaystyle \left\{(\tau_k)_1^{N},(w_{k})_1^{N}\right\}$ be a finite family of times and initial data 
of 
\eqref{2.1} associated to $u_0,$ and let $\displaystyle w:\big(0,1\big)\times\bigcup_{k=1}^N[S_k,T_k]
\longrightarrow\R$ be the 
solution of \eqref{2.1}. For every $\delta>0,$ there exists $\s_\delta=(\s_k)_1^N\in \R^N_+,\;v_\delta\in L^\infty((0,1)\times(0,T_N))$
such that, denoting by $u_\delta:(0,1)\times[0,T_N]\rightarrow\R$ the solution of \eqref{1.1} with bilinear control $v_\delta,$ we have
 \begin{equation}\label{Lem1}
\|u_{\delta}(\cdot,T_{k})-w(\cdot, T_{k})\|_{L^2(0,1)}
\leq\delta, \;\;\;	\forall k=1,\ldots,N.
\end{equation}
 \end{lem}
 {\bf Proof.} Fix $\displaystyle \left\{(\tau_k)_1^{N},(w_{k})_1^{N}\right\}$ and $\delta>0.$ 
 Let us consider the partition of $[0,T_N]$ in $2 N$ intervals introduced in 
 \eqref{N2}.
We will show that the bilinear control $v_\delta$ has the following expression
 $$
 v_\delta(x,t)=
  \begin{cases}
\quad   v_k^\delta(x,t) \;\;\;\;\;\;\,
 {\rm in} \;\;\;Q_{\mO_k
 } = (0,1) \times [T_{k-1},S_k] 
 , \;\:k=1,\ldots,N,
\\
\quad 0 
\qquad\quad\quad\;\;{\rm in} \,\;\;\;Q_{\mE_k
 } = (0,1) \times [S_k,T_k],
  \;\:\;\;\;\;k=1,\ldots,N. 
\end{cases}
$$
{\bf Step 1:} {\it
A useful energy estimate on $
(0,1)\times\displaystyle \bigcup_{k=1}^N [S_k,T_k]
\;\;\; (v_\delta\equiv0).$} \\
In the following, for every $k=1,\ldots,N,$ we will consider the following problem on $Q_{\mE_k}$
\begin{equation}\label{prk}
 \begin{cases}
\quad   u_t \; = \; u_{xx} \; 
 \; + \; f(u),
&\quad  {\rm in} \;\;\;Q_{\mE_k} = (0,1) \times [S_k, T_k], \;\; 
\\
\quad u (0,t) = u (1,t) = 0,
&\quad\quad\quad\quad\quad\quad\quad\;\; t \in [S_k, T_k],
\\
\quad u\:\mid_{t = S_k} \; = w_k+p_k, 
\end{cases}
\end{equation}
where $p_k\in H^{1}_0(0,1)$ are given functions,
and we will represent the 
solution 
of \eqref{prk}
as the sum of two functions $ w (x,t) $ and $h (x,t)$, which solve the following  problems in $Q_{\mE_k}$
\begin{equation}\label{lem h}
   \begin{cases}
\quad   w_t \; = \; w_{xx} \; 
\; + \; f(w)\;\;\;\;\text{ in }\;Q_{\mE_k},
\\
\quad w (0,t) = w (1,t) = 0,
\\
\quad w \:\mid_{t = S_k} \; = w_k\in C^{2+\vt} ([0,1] ), 
\end{cases}
\;
   \begin{cases}
\quad   h_{t} \; = \; h_{xx} \;
\; + \; (f(w+h) - f(w))\;\;\text{ in }\;Q_{\mE_k},
\\
\quad h (0,t) = h (1,t) = 0,
\\
\quad h \:\mid_{t = S_k} \; = p_k. 
\end{cases}\!\!\!\!\!
\end{equation}

{\it Evaluation of $ \|h(\cdot, T_k)\|_{L^2 (0,1)},\; k=1,\ldots,N$.} 
Let us fix $k=1,\ldots,N.$ Multiplying by $ h$ each member of the equation in the second problem of \eqref{lem h} and integrating by parts over $Q_{[S_k, t]} = (0,1) \times [S_k, t],\;\,t\in [S_k, T_k],$ 
keeping in mind (\ref{1.2a}), it follows that
\begin{align*}\frac{1}{2}\int_{ S_k}^{t}\int_0^1 (h^2)_t dx\,ds &=
\int_{ S_k}^{t}\int_0^1 h_{xx}h dx\,ds+\int_{ S_k}^t  \int_0^1 (f(w+h) -  f(w))h \,dx ds \\
&\leq-\int_{ S_k}^{t}\int_0^1 h^2_{x}dx\,ds+
L \int_{ S_k}^{t}  \int_0^1  h^2 dx ds,\quad \forall t\in [S_k, T_k].
\end{align*}
Then,
$$
\int_0^1 h^2 (x,t) dx 
\leq
 \int_0^1 p_k^2 (x) dx
+ 2L\int_{ S_k}^{t}  \int_0^1  h^2 dx ds, \quad \forall t\in [S_k, T_k].
$$
So, applying Gronwall's inequality we deduce
\begin{equation}\label{hGron}
\parallel h(\cdot,T_k) \parallel_{L^2 (0,1)} \; 
\leq  \; e^{L\widetilde{T}} \parallel  p_k \parallel_{L^2 (0,1)},\quad \text{ with }\;\; \displaystyle\widetilde{T}:=\sum_{k=1}^{N}\tau_{k}.
\end{equation}
{\bf Step 2:} {\it Steering.}
\begin{itemize}
\item {\it Steering the system from $u_0$ to $w_1
$ on $
[0,S_1]$.} 
  Applying Theorem \ref{th preserving} (with $k=1$ and $r_0=0$ in its statement), for every $\eta_1>0$ 
there exist $\s_1=\s_1(\eta_1, u_{0}, w_1)>0 
$
 and 
 a piecewise static bilinear control $v_1=v_1(\eta_1, u_{0}, w_1)\in L^\infty((0,1)\times(0,S_1))\;  (\text{with } 
 S_1=\s_1)$\; 
such that
\begin{equation}\label{ep i}
\|\overline{U}_1 (\cdot, S_1) - w_{1}(\cdot)\|_{L^2(0,1)}\leq
\eta_1,
\end{equation}
where 
 $\overline{U}_1$ is the solution of \eqref{u pres} on $(0,1)\times[0,S_1].$ 
Let us set 
\begin{equation}\label{pi1}
p_1(\cdot):=\overline{U}_1 (\cdot, S_1) - w_{1}(\cdot),
\end{equation}
we have 
$\overline{U}_1 (\cdot, S_1) = w_{1}(\cdot)+p_1(\cdot).$
\item {\it 
Steering the system from $w_1+p_1$ to $u_1
$ on $
[S_1,T_1]$.} 
 We consider the problem \eqref{prk}, written for $k=1.$ Due to \eqref{lem h}, we have 
 $u(\cdot,T_1)=w(\cdot,T_1)+h(\cdot,T_1),$
 where $u, w, h,$ defined on $\displaystyle(0,1)\times\bigcup_{k=1}^N[S_k,T_k],$ are the  solutions of the problem \eqref{prk}, the first problem in \eqref{lem h} and the second problem in \eqref{lem h}, respectively.
Thus, let us set 
\begin{equation}\label{re1}
u_1(\cdot):=w(\cdot,T_1)\qquad \text{ and } \qquad r_1(\cdot):=h(\cdot,T_1).
\end{equation}
\item{\it Steering the system from $u_{k-1}+r_{k-1}$
 to $w_k$ on $[T_{k-1},S_k],\; k=2,\ldots,N$.}\\
Similarly to 
\eqref{re1}, let us set
$$
u_{k-1}(\cdot):=w(\cdot,T_{k-1})\qquad \text{ and } \qquad r_{k-1}(\cdot):=h(\cdot,T_{k-1}).
$$
 Applying Theorem \ref{th preserving}, for every $\eta_k>0$ 
there exist $\s_k=\s_k(\eta_k, u_{k-1}, 
 w_k)>0$ 
 and 
 a piecewise static bilinear control $v_k=v_k(\eta_k, u_{k-1}, 
  w_k)\in L^\infty((0,1)\times(T_{k-1},S_k))\;  (\text{with } 
 S_k=T_{k-1}+\s_k),$\;
such that
\begin{equation}\label{ep k}
\|\overline{U}_k (\cdot, S_k) - w_{k}(\cdot)\|_{L^2(0,1)}\leq
\eta_k+ C_k
\|r_{k-1}\|_{L^2(0,1)},
\end{equation}
 where 
  $\overline{U}_k$ is the solution of \eqref{u pres} on $(0,1)\times[T_{k-1},S_k],$ 
 $C_k=C(u_{k-1}, w_k)\geq1$ is a constant, 
 as in Theorem \ref{th preserving}.
Moreover, we note that 
$
\overline{U}_k (\cdot, S_k) = w_{k}(\cdot)+p_k (\cdot),
$
where 
\begin{equation}\label{pik}
p_k(\cdot):=\overline{U}_k (\cdot, S_k) - w_{k}(\cdot).
\end{equation}
\item {\it 
Steering the system from $w_k+p_k$ to $u_k
$ on $
[S_k,T_k],\; k=2,\ldots,N$.} 
 Let us consider the problem \eqref{prk}.
 Due to \eqref{lem h}, we have 
 $u(\cdot,T_k)=w(\cdot,T_k)+h(\cdot,T_k),$
 where $u, w, h,$ defined on $\displaystyle(0,1)\times\bigcup_{k=1}^N[S_k,T_k],$ are the  solutions of the problem \eqref{prk}, the first problem in \eqref{lem h} and the second problem in \eqref{lem h}, respectively.
Thus, let us set 
\begin{equation}\label{rek}
u_k(\cdot):=w(\cdot,T_k)\qquad \text{ and } \qquad r_k(\cdot):=h(\cdot,T_k).
\end{equation}
\end{itemize}
 {\bf Step 3:} {\it Conclusions.} 
Given $\underline{\eta}=(\eta_1,\ldots,\eta_N)\in\R^N_+,$ let $\overline{U}_1,\ldots,\overline{U}_N$ be the solutions of \eqref{u pres} satisfying \eqref{ep i} and \eqref{ep k}.
Define 
 \begin{equation}\label{u delta}
 u_{\underline{\eta}}(x,t)=  \begin{cases}
\quad   \overline{U}_k(x,t) \;\;\;\;\;\;\,
\qquad\qquad\qquad\: {\rm in} \;\;\;  
 (0,1) \times [T_{k-1},S_k]
 , \;\:k=1,\ldots,N,\\
\quad w(x,t)+h(x,t)
\qquad\quad\quad\;\;{\rm in} \,\;\;\;
 (0,1) \times\displaystyle \bigcup_{k=1}^N[S_{k},T_k], 
\end{cases} 
(\footnote{We note that $\overline{U}_k, \;k=1,\ldots,N,$ depend on $\eta_k$ and $h$ depends on $\underline{\eta},$ but $w$ is independent of $\underline{\eta}.$
 }) 
\end{equation}
and observe that $u_{\underline{\eta}}$ is the solution of \eqref{1.1} corresponding to the piecewise static bilinear control\\
 $$
 v_{\underline{\eta}}(x,t)=
  \begin{cases}
\quad   v_k
(x,t) \;\;\;\;\;\;\,
 {\rm in} \;\;\;
 (0,1) \times [T_{k-1},S_k], 
  \;\:k=1,\ldots,N,
\\
\quad 0 
\qquad\quad\quad\;{\rm in} \,\;\;\;
\displaystyle (0,1) \times \bigcup_{k=1}^N[S_{k},T_k].
\end{cases}
$$
 Thus,
to show \eqref{Lem1} 
it is sufficient to prove by induction that there exists $\underline{\eta}(\delta)=(\eta_1(\delta),\ldots,\eta_N(\delta))\in\R^N_+,$ such that, for the corresponding $u_{\delta
}:=u_{\underline{\eta}(\delta)}, 
$ 
as in \eqref{u delta}, 
the following inequality holds
\begin{equation}\label{ind}
\|\overline{U}_k(\cdot,S_k)-w_k(\cdot)\|_{L^2(0,1)}\leq\frac{\delta\,C_k}{\small2^{N-k}e^{(N-k)L\widetilde{T}}{\small\small\displaystyle\prod_{h=k}^{N}}C_{h}}, \;\;\forall k=1,\ldots,N.
\end{equation}
Indeed, by \eqref{u delta}, \eqref{hGron}, \eqref{pi1} and \eqref{pik}, and \eqref{ind} we obtain \eqref{Lem1} and 
complete the proof.\\
{\bf Step 4:} {\it 
Proof of \eqref{ind}.} 
\begin{itemize}
\item{\it Base case.} 
Choosing $\eta_1=\eta_1(\delta):=\dfrac{\delta\,C_1}{2^{N-1}e^{(N-1)L\widetilde{T}}{\displaystyle\prod_{h=1}^{N}}C_{h}}
$ in \eqref{ep i}, we obtain \eqref{ind} for $k=1.$\\
\item
{\it Inductive step.} Let $k=1,\ldots,N-1$. By the inductive assumption, let us suppose that the inequality \eqref{ind} holds for the index $k.$ Thus, we will prove \eqref{ind} for the index $k+1.$ 
Choosing $\eta_{k+1}=\eta_{k+1}(\delta):=\dfrac{\delta\,C_{k+1}}{2^{N-k}e^{(N-k-1)L\widetilde{T}}{\displaystyle\prod_{h=k+1}^{N}}C_{h}}$ in \eqref{ep k}, keeping in mind \eqref{re1} or 
\eqref{rek}, \eqref{hGron}, \eqref{pi1} or \eqref{pik}, and the induction assumption we deduce
\begin{align*}
\|\overline{U}_{k+1}(\cdot,S_{k+1})-w_{k+1}(\cdot)\|_{L^2(0,1)}
&\leq \eta_{k+1}+C_{k+1}\|r_{k}\|_{L^2(0,1)}
\\
&\leq\eta_{k+1}+C_{k+1}e^{L\widetilde{T}}\|\overline{U}_{k}(\cdot,S_{k})-w_{k}(\cdot)\|_{L^2(0,1)}
\\
&\leq \eta_{k+1}+C_{k+1}e^{L\widetilde{T}}\frac{\delta\,C_k}{\small2^{N-k}e^{(N-k)L\widetilde{T}}{\small\small\displaystyle\prod_{h=k}^{N}}C_{h}}\\
&=\frac{\delta\,C_{k+1}}{2^{N-k-1}e^{(N-k-1)L\widetilde{T}}{\small\displaystyle\prod_{h=k+1}^{N}}C_{h}},
\end{align*}
from which the inequality \eqref{ind} is proved. $\;\qquad\qquad\qquad\qquad\qquad\qquad\qquad\qquad\diamond$\\
\end{itemize}

Now, we can prove Theorem \ref{th1.1}.\\
{\bf Proof (of Theorem \ref{th1.1}).}
Let us fix $ \eta > 0.$
Let $ u_{0}\in 
H^1_0(0,1)$ have $n$ points of sign change at $ x^0_{l} \in (0,1),\, 
$ 
with $0=:x^0_0<x^0_l<x^0_{l+1}\leq x^0_{n+1}:=1,\;\; l=1,\ldots,n$.
Let $ u^*\in H_0^1 (0,1)$ have $n$ points of sign change at $ x^*_{l} \in (0,1),\,  l=1,\ldots,n,$ such that $0:=x^*_0<x^*_l<x^*_{l+1}\leq x^*_{n+1}:=1.$ \\
{\bf Step 1:} Applying Theorem \ref{th2.2}, for every $\ve>0$ there exist $N_{\ve}\in\N$ and a finite family of times and initial data 
$ \left\{(\tau_k)^{N_{\ve}}_1, (w_{k})^{N_{\ve}}_1\right\}$ 
 such that, for any $(\s_k)^{N_\ve}_1\in\R^{N_\ve}_+,$   
the 
solution\\ $\displaystyle w^{\ve}:\big(0,1\big)\times\bigcup_{k=1}^{N_\ve}[S_k,T_k]\longrightarrow\R$ of 
problem \eqref{2.1} satisfies 
\begin{equation}\label{w eps}
 w^{\ve} (x, T_{N_{\ve}}
 )=0 \qquad \Longleftrightarrow \qquad x=x_l^{\ve},\;\;l=0,\ldots,n+1,  
 \end{equation}
for some points $0:=x^{\ve}_0<x^{\ve}_l<x^{\ve}_{l+1}\leq x^{\ve}_{n+1}:=1,\;l=1,\ldots,n,$ such that
\begin{equation}\label{x eps}
\displaystyle\sum_{l=1}^n|x^*_l-x^{\ve}_l|<\ve.
 \end{equation}
Moreover, $w^\ve(\cdot, T_{N_{\ve}})$ has the same order of sign change as for $u_0.$ \\
{\bf Step 2:} Applying Lemma \ref{pre proof th}, for every $\delta>0$ there exists $\s_{\ve,\delta}=(\s_k)_1^{N_\ve}\in \R^N_+,\;v_{\ve,\delta}\in L^\infty((0,1)\times(0,T_{N_\ve}))$
such that, denoted by $u_{\ve,\delta}:(0,1)\times[0,T_{N_\ve}]\longrightarrow\R$ the solution of \eqref{1.1} with bilinear control $v_{\ve,\delta},$ we have
\begin{equation}\label{u eps delta}
u_{\ve,\delta}(\cdot,T_{N_\ve})=w^{\ve} (\cdot, T_{N_{\ve}})+
r_{N_\ve}(\cdot),
\end{equation}
and
\begin{equation}\label{r eps delta}
\|r_{N_\ve}\|_{L^2(0,1)}
\leq 
\delta\,.
\end{equation}
{\bf Step 3:} {\it Steering the system from $
u_{\ve,\delta}(\cdot,T_{N_\ve})$
 to $u^*.$} 
By  
\eqref{x eps} it is easy to show that there exist $\ve^*=\ve^*(\eta)>0$ and  $u_{\ve}^*\in H^1_0(0,1)$ such that, for every $\ve\in(0,\ve^*),$ we have
 \begin{equation}\label{w*}
 u_{\ve}^*(x)=0\Longleftrightarrow x=x_l^\ve,\;\;l=0,\ldots,n+1\quad\quad\textit{ and }\quad\quad \|u^*
 -u_{\ve}^*\|_{L^2(0,1)}\leq\frac{\eta}{3}\,.
 \end{equation}
 Since 
 $w^{\ve} (\cdot, T_{N_{\ve}})$ and $u^*_\ve$ have the same points of sign change (see \eqref{w eps} and \eqref{w*}), keeping in mind \eqref{u eps delta} 
 we can steer the system \eqref{1.1} from $u_{\ve,\delta}(\cdot,T_{N_\ve})$ 
 to $u^*_\ve$ at some time $T>T_{N_\ve}$.
 Indeed applying Theorem \ref{th preserving}, for every $\overline{\eta}>0$ there exist $\s^*=\s^*(\overline{\eta}, \ve, u_0, u^*)
 >0 
$ 
 and 
a piecewise static bilinear control $v^*=v^*(\overline{\eta}, \ve, u_0, u^*)
 \in L^\infty((0,1)\times(T_{N_\ve},T)),$ with $T:=T_{N_\ve}+\s^*,$ 
such that, using also the inequality \eqref{r eps delta}, we obtain
\begin{equation}\label{final}
\|U^* (\cdot, T) - u^{*}_\ve(\cdot)\|_{L^2(0,1)}
\leq \overline{\eta}+ 
 C^*(\ve)\|r_{N_\ve}(\cdot)\|_{L^2(0,1)}\leq \overline{\eta}+ 
 C^*(\ve)\delta,
\end{equation}
where 
 $U^*$ is the solution of \eqref{u pres} on $(0,1)\times[T_{N_\ve},T]$ with initial state $u_{\ve,\delta}(\cdot,T_{N_\ve}),$ 
 $C^*(\ve)=C(\ve, u_0, u^*
 )
 .$\\
{\it Conclusions.} 
Let us fix $\ve:=\frac{\ve^*(\eta)}{2}$ and $\overline{\eta}:=\frac{\eta}{3}.$ Then, we consider the constant $C^*(\ve):=C^*\left(\frac{\ve^*(\eta)}{2}\right)$ of \eqref{final} and we choose $\delta:=\frac{\eta}{3C^*\left(\frac{\ve^*(\eta)}{2}\right)}.$  So, by \eqref{w*} and \eqref{final}, we obtain the conclusion
$$\|U^*(\cdot, T) - u^{*}(\cdot)\|_{L^2(0,1)}\leq\|U^* (\cdot, T) - u^{*}_\ve(\cdot)\|_{L^2(0,1)}+\|u^*- u^{*}_\ve\|_{L^2(0,1)} \leq 
  \eta.\qquad\qquad\diamond$$
  
 
\section{
Proof of Theorem \ref{th2.2} }\label{ADCP} 

{The plan of this section is as follows: 
}
\begin{description}
\item
In Section \ref{pre lem}, by Lemma \ref{l1} we construct suitable {\it initial data $w_k$'s}  to be used in the proof of Theorem \ref{th2.2}.
By Lemma \ref{p1} we construct the n {\it curves of sign change} associated to the 
the n initial points of sign change
and we prove some useful properties for the construction of suitable 
control strategies.
\item In Section \ref{OPS}, we construct a suitable particular family of {\it times and initial data},
that allows to move the $n$ initial points of sign change towards the $n$ target points of sign change. 
In this section, we also introduce
the definitions of {\it gap and target distance functional}, defined on the set of the {order processing steering times and initial data} 
at the beginning of this section). 
\item In Section \ref{Th3}, after proving a technical proposition (Proposition \ref{P2}) we show how to steer the points of sign change of the solution
arbitrarily close to the target points. 
\end{description}
For the notation of this section we refer to 
Section \ref{pure diffusion}.
\subsection{Preliminary results}\label{pre lem}
Let us prove the following lemma.
\begin{lem}[Construction of suitable initial data $w_k$'s]\label{l1} 
Let $x_l\in [0,1],\;l=0,\ldots,n+1,$ be such that 
$0=x_0<x_1<\cdots<x_n<x_{n+1}=1.$
Let $\alpha=(\alpha_0,\ldots, \alpha_{n+1})\in \R^{n+2},\;\beta=(\beta_0,\ldots, \beta_{n+1})\in \R^{n+2}\;$ be such that $\alpha_{l}\,\alpha_{l+1}<0, \alpha_l\in\{-1,1\},\;\beta_l\in\{-1, 0, 1\}, \,l=0,\ldots,n,\;  \beta_0  =  \beta_{n+1} =0$.
Let $\displaystyle \widetilde{\rho}=
\min_{l=0,\ldots,n}\left\{x_{l+1}-x_l
\right\}$. 
Then, there exists $w\in C^{^\infty}([0,1])$ such that
\begin{itemize}
\item[$\star$] $w(x)=0\;\;\Longleftrightarrow\;\; x=x_l,\; 
 l=0,\ldots, n+1;$
\item[$\star$] $w'(x_l)=\alpha_l,$\; $w''(x_l)=\beta_l,\;
l=0,\ldots, n+1;$
\item[$\star$] $\|w\|_{C^{k}([0,1])}\leq C(k,\widetilde{\rho}), \;\forall k\in \N$.
\end{itemize}
\end{lem}
{\bf Proof.}  For every $l=0,\ldots,n+1,$ set
$$
v_l(x)=\alpha_l(x-x_l)+\frac{\beta_l}{2}(x-x_l)^2,\;\forall x\in \R.$$
Note that  each $v_{l}(x)$ has no critical points in $\left[x_l -\frac{\widetilde{\rho}}{2}, x_l + \frac{\widetilde{\rho}}{2}\right]$.
Set $\rho=\frac{1}{2} \widetilde{\rho}$  and  define 
$$ w(x)=\sum_{l=0}^{n+1}
\eta_l^\rho (x)\,v_l(x)+\sum_{l=0}^{n}
\alpha_l\mathds{1}_{[x_l,x_{l+1})}(x)\left[1-(\eta_l^\rho (x)+\eta_{l+1}^\rho (x))\right],\;\;\;x\in[0,1],$$
where $ \mathds{1}_{A}(x)  $ is the characteristic function of a set $A,$ and 
$\eta_j^\rho\in C^\infty(\R),\;j=0,\ldots,n+1,$ are  such that $\eta_j^\rho(x)=\eta_0^\rho(x-x_j)$ and $\eta_0^\rho$ has the following properties:
\begin{itemize}
  \item $\eta_0^\rho (-x)=\eta_0^\rho (x), \, 0\leq \eta_0^\rho(x)\leq 1,\: 
  \forall x\in \R;\; $
 $\eta_0^\rho (x)=1,\, \forall x\in [0,\frac{\rho}{2}];$\; $\eta_0^\rho (x)=0,\, 
  \forall x\in[\rho,+\infty);$
  \item $\left|\dfrac{d^h\eta_0^\rho(x)}{dx^h} \right|\leq \frac{C_h}{\rho^h},\;\forall x\in \R$,  where $C_h$ is a positive constant and  $h\in \N.$
\end{itemize}
Observe that, for $x\in[x_l,x_{l+1}], \,l=0,\ldots,n,$
$$w(x)=\left\{ \begin{array}{ll}
v_l(x), \;\;\;\; &  {\rm if} \;\; x \in [x_l, x_l +\rho/2],   \\
\eta_0^\rho(x-x_l)v_l(x)+
\alpha_l[1-\eta_0^\rho(x-x_l)
], \;\;\;\; &  {\rm if} \;\;x \in (x_l   +  \rho/2, x_l +\rho), \\
\alpha_l,
&  {\rm if} \;\; x \in [x_{l} + \rho, x_{l+1}-\rho],   \\
\eta_0^\rho(x-x_{l+1})v_{l+1}(x)+
\alpha_l[1
-\eta_0^\rho(x-x_{l+1})], \;\;\;\; &  {\rm if} \;\; x \in (x_{l+1} -\rho, x_{l+1}-\rho/2),   \\
v_{l+1}(x), \;\;\;\; &  {\rm if} \;\;x \in (x_{l+1} - \rho/2, x_{l+1}]\,.
\end{array}
\right. 
$$ 
Notice that $w$ is of class $C^\infty(\R)$ by construction.
Moreover, our choice of  $\rho$ ensures that 
$w(x)$ has no points of sign change in $(x_{l}-\rho,x_{l})$ or in $(x_{l},x_{l}+\rho).$
This ends the proof of Lemma \ref{l1}.
$\qquad\qquad\qquad\qquad\qquad\qquad\qquad\diamond$

\setlength{\unitlength}{1.2mm}
\begin{picture}(150,55)(0,0)
\linethickness{1pt}
\put (2,25){\vector(1,0){105}} 

\put (10,4){\vector(0,1){45}}

\qbezier(10,25)(12,25.5)(13,27)
\qbezier(13,27)(17,39.2)(21,39.5)

\qbezier(25,39.5)(29,38)(31,27)

\qbezier(31,27)(32,25)(33.5,23.5)

\qbezier(33.2,23.5)(34,24)(36.5,11)

\qbezier(36.5,11)(37.5,8)(40.5,8)

\qbezier(47,8)(50,8)(51,12)

\qbezier(51,12)(54,23.5)(54.9,23.5)

\qbezier(54.9,23.5)(56.2,25)(59,27)

\qbezier(59,27)(61,29)(64.2,37)

\qbezier(64.2,37)(65.3,39.5)(69.3,39.5)

\qbezier(75.3,39.5)(79,39)(82,27.5)

\qbezier(82,27.5)(83,25.9)(84.2,25.3)

\put(30,20){\makebox(0,0)[b]{$x_{1}$}}
\put(86,28){\makebox(0,0)[b]{$1$}}

\put(6,39){\makebox(0,0)[b]{$1$}}
\put(6,8){\makebox(0,0)[b]{$-1$}}

\put(9,40){\line(1,0){1}}
\put(9,9){\line(1,0){1}}

\put(21,39.5){\line(1,0){4}}
\put(40,8){\line(1,0){8}}

\put(69.3,39.5){\line(1,0){6}}

\put(56.5,25){\line(0,1){1}}

\put(84.5,25){\line(0,1){1}}
\put(32.3,25){\line(0,1){1}}
\put(59,20){\makebox(0,0)[b]{$x_{2}$}}

\put(12,47){\makebox(0,0)[b]{$w$}}
\put(100,20){\makebox(0,0)[b]{$x$}}
\put(8,20){\makebox(0,0)[b]{$0$}}

\put(50,0){\makebox(11,10)[b]{{\rm Figure 1: 
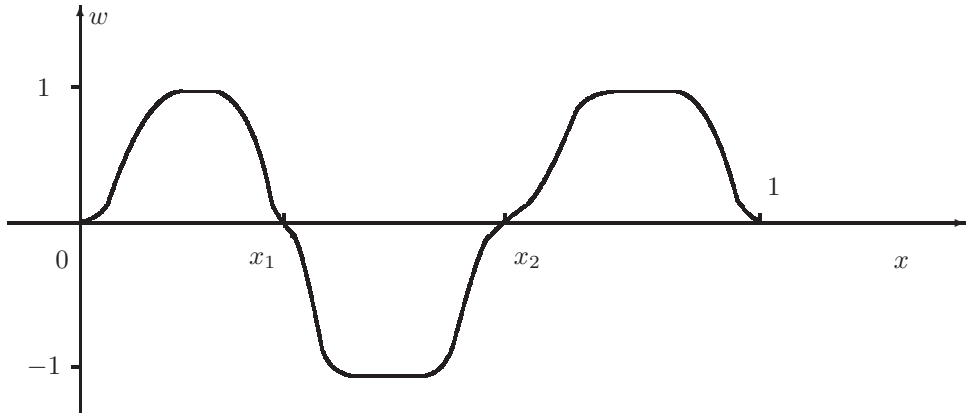
A function $w$ as in Lemma \ref{l1} with $2$ points of change of sign.}}}

\end{picture}

\begin{rem} \rm
In the above, we can costruct $\eta_0$ by the following expression on $\left(\frac{\rho}{2}, \rho\right)$
$$\eta_0(x)=\frac{e^{\frac{1}{(x-\rho)(x-\frac{\rho}{2})}}}{e^{\frac{1}{(x-\rho)(x-\frac{\rho}{2})}}
+e^{-\frac{1}{(x-\frac{\rho}{2})^2}}}=e^\frac{1}{1+h_0(x)},\;\; \mbox{ with }\; 
h_0(x)=\frac{-(2x-\frac{3}{2}\,\rho)}{(x-\frac{\rho}{2})^2(x-\rho)}.
$$
\end{rem}


In the whole section, 
$\vartheta\in (0,1)$ denotes the number that was fixed at the beginning of Section \ref{Sec Control Strategy}.
\begin{lem}[Construction of  the
curves of sign change]\label{p1}
Let $\alpha=(\alpha_0,\ldots, \alpha_{n+1})\in \R^{n+2}$ be 
such that $\alpha_{l}\,\alpha_{l+1}<0,\, \alpha_l\in\{-1,1\},
\,\alpha_{n+1}=-\alpha_{n},\,l=0,\ldots,n\,.\; $ 
Let 
$\widetilde{\rho}>0$ be. Let $x_l\in [0,1],\, l=0,\ldots,n+1,$ be such that 
$0=x_0<x_1<\cdots<x_n<x_{n+1}=1$ and $\displaystyle
\min_{l=0,\ldots, n}\left\{x_{l+1}
-x_l
\right\}=\widetilde{\rho}.$
Let $w_k\in C^{2+\vt}([0,1])$ be such that $w_k(x)=0$ if and only if
$x=x_l,\; 
 w_k'(x_l)=\alpha_l,
\,l=0,\ldots, n+1,$ $w_k''(0)=w_k''(1)= 0$ and
$\|w_k\|_{C^{2+\vt}([0,1])}\leq c
,$ for some positive constant $c=c(\widetilde{\rho}).$
Let $T>0$ and let $w\in C^{2+\vt,1+\frac{\vt}{2}}(\overline{Q}_T)$ be the solution of
the problem
\begin{equation}\label{u}
\begin{cases}
\quad   w_t \; = \; w_{xx} \;  + \; f(w) \;\;\;\;
 \qquad{\rm in} \;\;\quad Q_{T} = (0,1) \times (0, T)
\\
\quad w (0,t) = w (1,t) =0
\;\;\;\; \qquad\qquad\;\,\;\, \;\;t \in (0, T) \\
\quad w(x,0)
=w_k(x)\qquad\qquad\qquad\;\,\quad\;\;\;x\in(0,1) \,.
\end{cases}
(\footnote{For the existence, uniqueness and regularity of problem \eqref{u} see Remark \ref{solution of dtp} and Theorem 6.1 in \cite{Lad} (pp. 452-453).}) 
\end{equation}
Then, for every $\rho\in(0,\widetilde{\rho}]$ there exist $\widetilde{\tau}=\widetilde{\tau}({\rho}) >0$ and $M=M(\rho)>0$ such that, for each $l=1,\ldots,n,$ 
there exists a unique solution $\xi_l 
:[0,\widetilde{\tau}]\longrightarrow\R\, 
$  
of the initial-value problem
$$\;\;\;\displaystyle\begin{cases}
 \dot{\xi_l}
(t) \; = \;-\frac{ w_{xx}(\xi_l(t),t)}{w_x(\xi_l(t),t)},\;\;\;\; 
 t\in[0,\widetilde{\tau}],
\\
 \xi_l
(0)=x_l, 
\end{cases}
$$
that satisfies $w(\xi_l(t),t)=0,\;\forall t\in [0,\widetilde{\tau}],$ and 
$$\xi_l
\in C^{1+\frac{\vartheta}{2}}([0,\widetilde{\tau}]),\qquad \qquad \|\xi_l\|_{C^{1+\frac{\vartheta}{2}}([0,\widetilde{\tau}])}\leq M, \qquad 
\qquad \|\xi_l (\cdot) - x_l \|_{C([0,\widetilde{\tau}])} <\frac{\rho
}{2}.$$
\end{lem}
\begin{rem}\label{separeted curves}\rm
In Lemma \ref{p1}, since $\|\xi_l (\cdot) - x_l \|_{C([0,\widetilde{\tau}])} < \dfrac{\rho}{2}$ for each $l=1,\ldots,n,$ we also have that 
\begin{equation}\label{max1}
0:=\xi_{0}(t)<\xi_{l}(t)<\xi_{l+1}(t)<\xi_{n+1}(t):=1,\;\;\;\forall t\in[0,\widetilde{\tau}], \;\forall l=1,\ldots,n-1
.
\end{equation}
\end{rem}
\begin{defn}\label{scgs}
We call the functions $\xi_l:[0,\widetilde{\tau}]\longrightarrow \R,\; l=1,\ldots,n,$ given by Lemma \ref{p1}, Curves of Sign Change associated to the set of initial points of sign change $X
=(x_1,\ldots,x_n)$.
\end{defn}
{\bf Proof (of Lemma \ref{p1}).}\; Let us fix $\rho\in(0,\widetilde{\rho}].$\\
{\bf Step 1:} {\it Uniform estimate for $w$.} 
Due to Theorem 6.1 of \cite{Lad} (pp. 452-453), the solution $w$ of \eqref{u} belongs to $C^{2+\vartheta,1+\frac{\vartheta}{2}}(\overline{Q}_T)$ (\footnote{One can note that the initial datum satisfies the compatibility condition $w_k''(0)=w_k''(1)=0,$ as in Theorem 6.1 of \cite{Lad}.})
and, for some constant $K=K( \| w_k \|_{C^{2+\vartheta}([0,1])})>0,$ depending only on $\| w_k \|_{C^{2+\vartheta}([0,1])}$ (see (6.8)-(6.12) on pp. 451-452 in \cite{Lad}), 
we have
$
\| w \|_{C^{2+\vartheta,1+\frac{\vartheta}{2}}(\overline{Q}_T)} \; \leq
K.$
Thus, since $ \| w_k \|_{C^{2+\vartheta}([0,1])}\leq  c(\widetilde{\rho}),$ 
we deduce that
\begin{equation}\label{uest}
\| w \|_{C^{2+\vartheta,1+\frac{\vartheta}{2}}(\overline{Q}_T)} \; \leq
K( \| w_k \|_{C^{2+\vartheta}([0,1])})
\leq \; C, 
\end{equation}
for some positive constant $C=C(\widetilde{\rho})$ depending only on $\widetilde{\rho}.$\\
{\bf Step 2:} {\it Existence and regularity of curves of sign change.}
 For any fixed $l=1,\ldots,n,$
 since $w_x(x_l,0)=\alpha_l\neq0$ and $w_x(x,t)$ is a continuous function in $(x_l,0)\in\overline{Q}_T,$  there exist $\delta_l\in\Big(0,\min\big\{\frac{1}{2C},\rho\big\}\Big)$(\footnote{$C$ is the constant present in \eqref{uest}.}) and
  $T_l>0$ 
 such that $w_x(x,t
)\neq0,\;\forall (x,t)\in [x_l-\delta_l, x_l+\delta_l]\times[0,T_l].$\\
Let 
$\displaystyle \delta:=\min_{l=1,\ldots,n}\delta_l$ be.
For every $l=1,\ldots,n,$
we consider the Cauchy problems 
\begin{equation}\label{ODE}
\begin{cases}
\quad   \dot{\xi_l}
(t) \; = \;-\frac{ w_{t}(\xi_l
(t),t)}{w_x(\xi_l
(t),t)},\;\;\;\;
 t>0,
\\
\quad \xi_l
(0)=x_l\;.
\end{cases}
\end{equation}
Observe that 
$F(x,t):=-\frac{ w_{t}(x,t
)}{w_x(x,t
)}
$ is continuous on $[x_l-\delta, x_l+\delta]\times [0,T_l].$ 
Therefore, for every $l=1,\ldots,n,$ the problem \eqref{ODE} has a solution $\xi_l$ of class 
$C^1$ 
 on some interval $[0,\tau_l],$ with $0<\tau_l\leq T_l
.$ 
Moreover, since
$F
\in C^{\frac{\vartheta}{2}}([x_l-\delta, x_l+\delta]\times [0,T_l]),$ we conclude that $\xi_l\in C^{1+\frac{\vt}{2}}([0,\tau_l]).$
Furthermore, 
$w(\xi_l(t),t)=0,\;\;\forall t\in[0,\tau_l],$ 
because 
$$\frac{d}{dt}w(\xi_l
(t),t)=w_{t}(\xi_l
(t),t)+w_{x}(\xi_l
(t),t)\dot{\xi_l}
(t)=0, \forall t\in [0,\tau_l],
\quad
\text{ and }
\quad
w(\xi
_l(0),0)=
w_k(x_l)=0.$$
Moreover, since
 $f(w(\xi_l(t)),t)=0,\;\forall t\in[0,\tau_l],$
we also have that
$$\displaystyle
\dot{\xi_l}
(t)=-\frac{ w_{t}(\xi_l
(t),t)}{w_x(\xi_l
(t),t)}
=-\frac{ w_{xx}(\xi_l
(t),t)}{w_x(\xi_l
(t),t)},\;\;\,\forall t\in[0,\tau_l].$$
{\bf Step 3:} {\it Uniform estimates for the curves of sign change.}
For any fixed $l=1,\ldots,n,$ we consider
the number $\delta=\delta(\rho)>0,$ $\displaystyle\delta=\min_{l=1,\ldots,n}\delta_l<\min\Big\{\frac{1}{2C},\rho\Big\},$  introduced in Step 2, 
and the uniform time $\widetilde{\tau}=\widetilde{\tau}(\rho)>0,$ 
 $$\widetilde{\tau}
 := \min\Big\{\Big(\frac{1}{2C}-\delta\Big)^\frac{2}{\vt}, \frac{\delta^2}{3}, \displaystyle \min_{l=1,\ldots,n} \tau_l\Big\}.\;\;(\footnote{We note that $\widetilde{\tau}=\widetilde{\tau}(\rho)$ is not dependent on $l.$})$$ 
We remember that the function $t\mapsto w_x(x,t)$ belongs to $C^{\frac{\vt}{2}}([0,\widetilde{\tau}]),$ and the function $x\mapsto w_x(x,t)$ belongs to $C^{1+\vt}([0,\widetilde{\tau}]).$ Thus, for every $(x,t)\in (x_l-\delta,x_l+\delta)\backslash\{x_l\}\times (0,\widetilde{\tau}),$ by \eqref{uest} we have
\begin{align}\label{uxest}
| w_{x}(x,t) - \alpha_l | &=  | w_{x}(x,t) - w_{x}(x_l,0) | 
\leq | w_{x}(x,t) - w_{x}(x,0)|  + | w_{x}(x,0) - w_{x}(x_l,0)|\nonumber\\
&= \frac{| w_{x}(x,t) - w_{x}(x,0)|}{t^{\frac{\vt}{2}}}\,t^{\frac{\vt}{2}}+\frac{| w_{x}(x,0) - w_{x}(x_l,0)|}{ |x - x_l|} |x - x_l|\nonumber\\
&\leq \|w\|_{C^{2+\vt,1+\frac{\vt}{2}}(\overline{Q}_T)}( t^{\frac{\vartheta}{2}}+|x - x_l|)\leq \;  
C ( t^{\frac{\vartheta}{2}}+|x - x_l| )
 \leq C ( \widetilde{\tau}^{\frac{\vartheta}{2}}+\delta ).
\end{align}
Since
 $\delta<\frac{1}{2C}$ and
 $\widetilde{\tau}\leq \big(\frac{1}{2C}-\delta\big)^\frac{2}{\vt},$
we have
$C (\widetilde{\tau}^{\frac{\vartheta}{2}} + \delta) \leq \frac{1}{2},
$
so by \eqref{uxest} we deduce
$$\Big||w_{x}(x,t)| - |\alpha_l|\Big|\leq| w_{x}(x,t) - \alpha_l |\leq C (\widetilde{\tau}^{\frac{\vartheta}{2}} + \delta) \leq\frac{1}{2}.$$
Therefore, for every $l = 1, \ldots, n,$ having in mind that $ | \alpha_l | = 1,$ we obtain
\begin{equation}\label{0.5}
|w_{x}(x,t)|\geq | \alpha_l |-\frac{1}{2}=\frac{1}{2}, \;\;\;\;\;\;\forall (x,t)\in (x_l-\delta,x_l+\delta)\times (0,\widetilde{\tau})\,.
\end{equation}
Then, by \eqref{uest} and \eqref{0.5}, keeping in mind that $ \widetilde{\tau}\leq\frac{\delta^2}{3}$ and $\delta<\min\{\frac{1}{2C},\rho\},$ we deduce 
\begin{multline}\label{delta1}
| \xi_l (t) - x_l |=\left| \int_0^t\dot{\xi_l} (s)\,ds \right|   \leq  \int_0^{\widetilde{\tau}} \frac{ | w_{xx}(\xi_l(s),s) | }{|w_x(\xi_l(s),s)| } ds \; 
\leq  \;\frac{\widetilde{\tau}  \|w\|_{C^{2+\vt,1+\frac{\vt}{2}}(\overline{Q}_T)}
 }{\displaystyle\min_{s \in [0, \widetilde{\tau}]}|w_x(\xi_l(s),s)|}
\\ \leq  \frac{\widetilde{\tau}  C}{\displaystyle\min_{s \in [0, \widetilde{\tau}]}|w_x(\xi_l(s),s)|}\leq\,2\widetilde{\tau} C<\frac{\widetilde{\tau}}{\delta}\leq\frac{1}{\delta}\frac{\delta^2}{3}<  \frac{\rho}{3} ,\;\;\;\forall t\in [0,\widetilde{\tau}].
\end{multline}  
{\bf Step 4:} {\it Uniqueness of the curves of sign change.}
We note that, although one cannot claim uniqueness for the Cauchy problem \eqref{ODE}, 
a posteriori the $\xi_l$'s 
 turn out to be uniquely determined.
 Indeed, setting $\xi_0(t)\equiv0,\;\xi_{n+1}(t)\equiv0,\;\forall t\in [0,\widetilde{\tau}],$ 
 one can apply the maximum principle for semilinear parabolic equations (see Remark \ref{r1.2}) on the domains 
 $$\left\{(x,t)|x\in\left[\xi_l(t), \xi_{l+1}(t)\right],\,t\in[0,\widetilde{\tau})\right\},$$ 
 for every $l=0,\ldots,n.$ The fact that the initial datum $w_k(x)$ doesn't change sign on $(x_l,x_{l+1})$ and the boundary conditions in 
  \eqref{u}
   imply that, for every $t^*\in [0,\widetilde{\tau}),$
$$w(x,t^*)=0\Longleftrightarrow \,x=\xi_l(t^*),\; l=0,\ldots,n+1,$$
this completes the proof of  Lemma \ref{p1}. 
$\qquad\qquad\qquad\qquad\qquad\qquad\qquad\qquad\qquad\qquad\qquad\diamond$\\
\begin{rem} \rm
An alternative proof of 
Lemma \ref{p1} can be obtained by using the implicit function theorem 
instead of solving problem \eqref{ODE}. 
\end{rem}
\subsection{Construction of Order Processing Steering sets of Times and Initial Data  
}\label{OPS}
In this section we define the set of {\it Order Processing Steering sets of Times and Initial Data}
that permit to move the points of sign change towards the desired targets. 
\subsubsection*{Notation}
Let us consider the initial state $u_0\in H^1_0(0,1).$ For simplicity of notation let us set $x^0_0:=0$ and $x^0_{n+1}:=1,$ and let us consider the set of $n$ points of sign change of $u_0,$
 $X^0=(x^0_1,\ldots,x^0_n),$ where
$0=x^0_0<x^0_l<x^0_{l+1}\leq  x^0_{n+1}=1,\;\; l=1,\ldots,n.$ Let $\displaystyle \rho_0=
\min_{l=0,\ldots,n}\big\{x^0_{l+1}-x^0_l
\big\}$. 
Let us define
\begin{equation}\label{lambda}
\lambda(x^0_l)=
\begin{cases}
 1,\;\;\;\;
 \text{ if }  u_0(x)>0 \text{ on } (x^0_l,x^0_{l+1}),
\\
-1,\;\, \text{ if } u_0(x)<0 \text{ on } (x^0_l,x^0_{l+1}),
\end{cases}
\!\!\!
 l=0,\ldots,n\;\; \mbox{ and } \;\;\lambda(x^0_{n+1})=-\lambda(x^0_{n}).\,\;(\footnote{Since $x_l^0,\,l=1,\ldots,n,$ are points of sign change, we note that $\lambda(x^0_{l+1})=-\lambda(x^0_{l}),\,l=0,\ldots,n.$})
\end{equation} 
For simplicity of notation let us set $x^*_0:=0$ and $x^*_{n+1}:=1,$ and let us also consider the set of $n$ target points 
$X^*=(x^*_1,\ldots,x^*_n),$ where
$0=x^*_0<x^*_l<x^*_{l+1}\leq x^*_{n+1}=1,\;\; l=1,\ldots,n.$ 
\subsubsection*{Order Processing Steering sets of Times and Initial Data} 
Let us consider the same $\vartheta\in (0,1)$ that was fixed at the beginning of Section \ref{Sec Control Strategy}.\\ 
Now, we define 
the set of the {\it Order Processing Steering Times and Initial Data} 
associated to $u_0$ and $X^*$ 
 $${\cal{W^*}}(u_0):=\left\{W^N\,|\;W^N=\big\{(\tau_k)^N_{1},(w_k)^N_{1}\big\}, N\in \N
 \right\},$$
 where a generic set 
 $\displaystyle W^N=\big\{(\tau_k)_1^N, (w_k)_1^N\big\}
 \in{\cal{W^*}}(u_0)$ 
  is defined below.
Let us fix any $(\s_k)_1^N\in\R^N_+$ and, recalling the notation introduced in 
\eqref{N1}-\eqref{N2} in Section \ref{Sec Control Strategy}, for $k=1,\ldots,N,$ 
we consider the initial value 
pure diffusion problems on disjoint time intervals 
\eqref{2.1} 
$$   \begin{cases}
\quad   w_t \; = \; w_{xx} \;  + \; f(w) \;\;\;\;
 {\rm in} \;\;\;Q_{\mE_k
 } = (0,1) \times [S_{k}, T_{k}
 ]
\\
\quad w (0,t) = w (1,t) = 0,
\qquad\qquad\qquad\;\;\:\, t \in [S_k, T_{k}]
\\
\quad w\mid_{t = S_k} \; = w_{k}(x),
\;\;\;\;\qquad w_{k}'' (x) \mid_{x = 0, 1} = 0. 
\end{cases}
$$
 \begin{description}
 \item \underline{\it Properties of $\big\{\tau_1, w_1\big\}$}.
By Lemma \ref{l1}, there exists $w_1\in C^{2+\vt}([0,1]),$ with $\|w_1\|_{C^{2+\vartheta}([0,1])}\leq c_1,$
 for some positive constant $c_1= c(
{\rho_0}
),$
 such that
\begin{itemize}
\item $w_1(x)=0\Longleftrightarrow x=x^0_l,\; 
 l=0,\ldots, n+1;$
\item $w_1'(x^0_l)=\lambda(x^0_l),\, 
w_1''(x^0_l)=-\lambda(x^0_l)\mu_1(x^*_l-x^0_l),\;l=0,\ldots, n+1,
$\\
where
$
\mu_1
(x_l^*-x_l^0):=sgn(x_l^*-x_l^0)=
\begin{cases}
 1,\;\;\;\;
 \text{ if } \,x_l^0<x_l^*,\;\;\;\;\\
0,\;\;\;\;\text{ if }\,x_l^0=x_l^*,\\
-1,\,\;\text{ if }\,x_l^0>x_l^* \,.
\end{cases}
$
\end{itemize}
Let $w
$ be the solution to
 $$
 \begin{cases}
\quad   w_t \; = \; w_{xx} \;  + \; f(w) \;\;\;\;
\qquad 
(x,t)\in(0,1) \times (S_1, +\infty)
\\
\quad w (0,t) = w (1,t) = 0
\; \quad\qquad\quad\; t \in (S_1, +\infty)
\\
\quad w(x,S_1)
 \, = w_1(x)\;\;\qquad\qquad\;\;\;\;\;\; x\in(0,1),\\ 
\end{cases}
$$
where $S_1=\s_1.$
By Lemma \ref{p1}, for every $\rho\in\left(0,\rho_0\right]$ there exist $\widetilde{\tau}_1=\widetilde{\tau}_1({\rho})>0,$ $M_1=M_1(
{\rho
})>0$ 
and
$n$ {\it small curves of sign change} (associated 
to the set of n initial points of sign change $X^{0}=(x^{0}_1,\ldots,x^{0}_n)$)
  $\xi^1_l
\in C^{1+\frac{\vartheta}{2}}([S_1,\widetilde{T}_1]), 
\,\widetilde{T}_1=S_1+\widetilde{\tau}_1,\, l=1,\ldots,n,$
such that $w(\xi^1_l
(t),t)= 0, \forall t\in [S_1,\widetilde{T}_1].
$ Moreover, 
\begin{equation}\label{xi}
\!\!\!\!\!
\begin{cases}
\quad  \dot{\xi^1_l}
(t) \; = \;-\frac{ w_{xx}(\xi^1_l
(t),t)}{w_x(\xi^1_l
(t),t)},
\;\;\;\;
 t\in[S_1,\widetilde{T}_1],
\\
\quad \xi_l^1
(S_1)=x^0_l\;,
\end{cases}
\mbox{ and   }\;\,
\|\xi^1_l
\|_{C^{1+\frac{\vartheta}{2}}([S_1,\widetilde{T}_1])}\leq M_{1}.
\end{equation}
Let us set
$\xi^1_0(t)\equiv 
0$ and $\xi^1_{n+1}(t)\equiv 
1$ on $[S_1,\widetilde{T}_1].$ Furthermore, for every $l=1,\ldots,n-1,$ by Remark \ref{separeted curves} we have
\begin{equation}\label{sccs1}
0=\xi^1
_0(t)<\xi^1
_l(t)<\xi^1
_{l+1}(t)<\xi^1
_{n+1}(t)=1,\; \forall t\in [S_1,\widetilde{T}_1].
\end{equation}
Let us introduce the {\it Inactive Set}
$$L^{0}_{IS}:=\{l\,|\,l\in\{1,\ldots,n\},\,x_l^0=x_l^*\}$$
and let us consider the set of the {\it stopping times}  
$$\Theta_{1}
:=\{s\in(0,\widetilde{\tau}_{1}]\,|\,\xi_{l}^1(S_1+s)=x_l^*,\;\text{ for some }\,l\in\{1,\ldots,n\}\backslash L^{0}_{IS}\}.
$$
Let us set
\begin{equation}\label{tau***1}
\displaystyle
\tau_{1}=\begin{cases}
\; \widetilde{\tau}_1
\qquad\qquad\quad \text{ if }\quad 
\;\;\Theta_{1}=\varnothing,
\\
\;\displaystyle\min \Theta_{1},
\qquad\qquad\text{ otherwise },
\end{cases}
\end{equation}
by \eqref{N1} we have
$T_{1}=S_1+\tau_1.$
 \item \underline{\it Properties of $\big\{\tau_k, w_k\big\},\;\;k=2,\ldots,N$.}
 By iterate application of Lemma \ref{l1} and Lemma \ref{p1}, one can set
$x^{k-1}_l:=\xi_l^{k-1}(T_{k-1}),\,l=1,\ldots,n$ and $X^{k-1}=(x_1^{k-1},\ldots,x_n^{k-1}),$ 
where $\xi_l^{k-1}(t), \;t\in[S_{k-1},T_{k-1}],$ are the $n$ {\it curves of sign change} associated to the initial state $w_{k-1}$ and to the set of n points of sign change $X^{k-2}=(x^{k-2}_1,\ldots,x^{k-2}_n)$. For simplicity of notation, let us set $x^{k-1}_0:=0$ and $x^{k-1}_{n+1}:=1.$ 
Let $\displaystyle \rho_{k-1}
=
\min_{l=0,\ldots,n}\left\{x^{k-1}_{l+1}-x^{k-1}_l
\right\}$.
\\
Let us introduce the 
{\it Inactive Set } 
 $$L^{k-1}_{IS}
:=\{l,\;l\in\{1,\ldots,n\}|\exists h_l\in\{1,\ldots,k-1\}:
x^{h_l}_l=x^*_l\},\; (\footnote{We note that $L^{k-1}_{IS}\subseteq\{1,\ldots,n\}$ 
 is a family of sets increasing in $k.$}) 
$$
  that is, the set of the indexes of the points of sign change that have already reached the corresponding target points of the sign change set $X^*=(x_1^*,\ldots,x_n^*)$ in some previous time instant. 
Then, let us set\\
$$
\displaystyle
\mu_k
(x_l^*-x_l^0)=
\begin{cases}
 0,\;\;\;\;
\;\;\quad\qquad\qquad\qquad \text{ if } \,l\in L_{IS}^{k-1}
,
 \;\;\;\;\\
 sgn(x_l^*-x_l^0)
\;\;\qquad\;\;\quad\text{ if }\, l\not\in L_{IS}^{k-1}\,.
\end{cases}
(\footnote{We observe that the definition of $\mu_k, \,k=2,\ldots,N,$ is consistent with the one of $\mu_1.$})
$$

Thus, by Lemma \ref{l1} we can choose $w_k\in C^{2+\vt}([0,1]),$ with $\|w_k\|_{C^{2+\vartheta}([0,1])}\leq c_k,$
 for some positive constant $c_k= c(
{\rho_{k-1}}
),$ 
such that
\begin{itemize}
\item $w_k(x)=0\Longleftrightarrow x=x^{k-1}_l,\; 
 l=0,\ldots, n+1;$ 
\item $w_k'(x^{k-1}_l)=\lambda(x^{0}_l),$\,
 \;$w_{k}''(x^{k-1}_l)=-\lambda(x^0_l)\mu_{k}(x^*_l-x^0_l),\;
l=0,\ldots, n+1.$ 
\end{itemize}

Let $w
$ be the solution to
 $$
 \begin{cases}
\quad   w_t \; = \; w_{xx} \;  + \; f(w) \;\;\;\;
\qquad 
(x,t)\in(0,1) \times (S_k, +\infty)
\\
\quad w (0,t) = w (1,t) = 0
 \quad\quad\qquad\;\; t \in (S_k, +\infty)
\\
\quad w(x,S_k)
 \, = w_k(x)\;\;\qquad\qquad\;\;\; \,x\in(0,1), 
\end{cases}
$$
where $\displaystyle S_k=\sum_{h=1}^{k-1}\left(\s_h+\tau_h\right)+\s_k.$
By Lemma \ref{p1}, for every $\rho\in(0,\rho_{k-1}]$ there exist $\displaystyle\widetilde{\tau}_k=\widetilde{\tau}_k(\rho)>0,$ $M_k=M_k(\rho)>0$ 
and $n$ {\it small curves of sign change} (associated 
to the set of $n$ intermediate points of sign change
$X^{k-1}$), $\xi^k_l
\in C^{1+\frac{\vartheta}{2}}([S_k,T_{k}]),\,\displaystyle 
\widetilde{T}_k=S_k+\widetilde{\tau}_k,\; l=1,\ldots,n,$
such that $ w
(\xi^k_l
(t),t)=0,\,\forall t\in [S_k,\widetilde{T}_k].$ Moreover, 
\begin{equation}\label{xik}
\!\!\!\!\!\begin{cases}
\quad   \dot{\xi^k_l}
(t) \; = \;-\frac{ w_{xx}(\xi^k_l
(t),t)}{w_x(\xi^k_l
(t),t)},
\;\;\;\;
 t\in[S_k,\widetilde{T}_{k}],
\\
\quad \xi_l^k
(S_k)=x^{k-1}_l\;,
\end{cases}
\mbox{ and   }\;\,
\|\xi^k_l
\|_{C^{1+\frac{\vartheta}{2}}([S_k,\widetilde{T}_{k}])}\leq M_{k}.
\end{equation}
Let us set
$\xi^k_0(t)\equiv 
0$ and $\xi^k_{n+1}(t)\equiv 
1$ on $[S_k,\widetilde{T}_{k}].$  Furthermore, for every $l=1,\ldots,n-1,$ by Remark \ref{separeted curves} we have
 \begin{equation}\label{sccsk}
 0=\xi^k
_0(t)<\xi^k
_l(t)<\xi^k
_{l+1}(t)<\xi^k
_{n+1}(t)=1,\; \forall t\in [S_k,\widetilde{T}_k].
\end{equation}
Let us introduce the set of the {\it stopping times} 
$$\Theta_{k}:=\{s\in(0,\widetilde{\tau}_{k}]\,|\,\xi_{l}^k(S_k+s)=x_l^*,\;\text{ for some }\,l\in\{1,\ldots,n\}\backslash L^{k-1}_{IS}\},$$
and let us set 
\begin{equation}\label{tau***k}
\tau_k=\begin{cases}
\;\widetilde{\tau}_{k}
\qquad\qquad \text{ if }\quad 
\;\;\Theta_{k}=\varnothing,
\\
\;\displaystyle\min \Theta_{k},
\qquad\qquad\text{ otherwise }.
\end{cases}
\end{equation}
Then, by \eqref{N1} we have $T_k=S_k+\tau_k.$
\end{description}
Now, we give some remarks about the above introduced set of the {order processing steering times and initial data} 
associated to $u_0$ and $X^*.$
\begin{rem}\label{remaining}\rm
We note that $\tau_k<\widetilde{\tau}_{k}$ for at most $n$ values of $k\in\{1,\ldots,N\}.$
\end{rem}
\begin{rem}\label{rem steer}
\rm
We note that, for each index $l\not\in L_{IS}^{k-1}\,(k=1,\ldots,N),$ 
by \eqref{xik} and the choice of the initial data $w_k$ we deduce
  $$ \dot{\xi^k_l}
(S_k) \; = \;-\frac{ w_{xx}(\xi^k_l
(S_k),S_k)}{w_x(\xi^k_l
(S_k),S_k)}=-\frac{-\lambda(x^0_l)\mu_k(x^*_l-x^0_l)}{\lambda(x^0_l)}=\mu_k(x^*_l-x^0_l),\;\;\;\;\xi_l^k
(S_k)=x^{k-1}_l.$$
Then, if $x^0_l<x^*_l$ we have $\dot{\xi^k_l}
(S_k) \; =\mu_k(x^*_l-x^0_l)>0,$ so the initial conditions $w_k$ permits to move the points of sign change $x^{k-1}_l$ on the right toward the desired $x^*_l.$  Otherwise, if $x^*_l<x^0_l$ the initial conditions $w_k$ permits to move the points of sign change on the left. 
\end{rem}
\begin{rem}\rm
We note that, for each inactive index $l\in L_{IS}^{k-1} \,(k=1,\ldots,N),$ 
we choose the initial data 
  such that the second derivative is equal to $0$ in the intermediate points of sign change $x^{k-1}_l.$ So, as 
 we will see later, such points will remain forever near the corresponding target points of sign change already reached  (we will see it in the next inequality \eqref{contr} in the proof of Proposition \ref{P2}). 
\end{rem}
\subsubsection*{Curves of Sign Change, Gap and Target Distance functional}
Given $W^N=\displaystyle \left\{(\tau_k)_1^{N},(w_{k})_1^{N}\right\}
\in{\cal{W^*}}(u_0),$ we introduce the {\it n Curves of Sign Change} 
associated to $W^N,$
as the functions $\displaystyle \xi_l^{W}:\bigcup_{k=1}^N
[S_k,T_k]
\longrightarrow\R, \;l=1,\ldots,n,$ such that 
$$\xi_l^{W}(t)=
\xi^{k
}_l
(t
),\qquad S_{k}\leq t\leq T_k,\;\; k=1,\ldots,N,$$
where any {curve}
$\xi^{k
}_l
,\;k=1,\ldots,N,$ is given by Lemma \ref{p1}  and is defined on $[S_k,T_k].$ 
For simplicity of notation, let us set 
$\xi^{W}_0(t)\equiv 
0$ and $\xi^{W}_{n+1}(t)\equiv 
1.$
Moreover, for every $l=1,\ldots,n-1,$ by \eqref{sccs1} and \eqref{sccsk},
we deduce that
$$0=\xi_0^{W}(t)<\xi_l^{W}(t)<\xi_{l+1}^{W}(t)< \xi_{n+1}^{W}(t)=1,\;\;\forall t\in  \displaystyle {\bigcup_{k=1}^N [S_k,T_k]
}
.$$
\begin{defn}
 We define the \textbf{gap functional} $\rho:{\cal{W^*}}(u_0)\rightarrow [0,1]$ in the following way
$$\rho(W^N)=\min_{l=0,\ldots,n}\min_{t\in
 \displaystyle {\mE}}\left\{\xi_{l+1}^{W}(t)-\xi_l^{W}(t)\right\},\quad \forall\, W^N=
 \left\{(\tau_k)_1^{N},(w_{k})_1^{N}\right\}
\in{\cal{W^*}}(u_0),$$
where $\displaystyle\mE:=\bigcup_{k=1}^N[S_k,T_k]
.$
\end{defn}
\begin{defn}
We define the \textbf{target distance 
functional} associated to the set $X^*,$\\ $J^*:{\cal{W^*}}(u_0)\rightarrow [0,1]$ such that
$$J^*(W^N)=
\sum_{l=1}^n \left|\xi_l^{W}(T_N)-x_l^*\right|,\,\quad \forall\, W^N=\displaystyle \left\{(\tau_k)_1^{N},(w_{k})_1^{N}\right\}\in{\cal{W^*}}(u_0).$$
\end{defn}
\subsection{End of the proof of Theorem \ref{th2.2}}\label{Th3}
\subsubsection*{Notation}
Let us consider the initial state $u_0\in H^1_0(0,1).$ 
Let us set $x^0_0:=0$ and $x^0_{n+1}:=1,$ and let us consider the set of $n$ points of sign change of $u_0,$
 $X^0=(x^0_1,\ldots,x^0_n)$ where
$0=x^0_0<x^0_l<x^0_{l+1}\leq  x^0_{n+1}=1,\;\; l=1,\ldots,n.$ Let $\displaystyle\rho_0=\min_{l=0,\ldots,n}\left\{x_{l+1}^0-x_{l}^0\right\}$ be.\\
Let us set $x^*_0:=0$ and $x^*_{n+1}:=1,$ and
let us consider the set of $n$ target points 
$X^*=(x^*_1,\ldots,x^*_n),$ where
$0=x^*_0<x^*_l<x^*_{l+1}\leq x^*_{n+1}=1,\;\; l=1,\ldots,n.$ \\
Let $\displaystyle\rho^*_0=\min_{l=0,\ldots,n}\left\{x_{l+1}^*-x_{l}^*,\;x_{l+1}^0-x_{l}^0\right\}$ be. Let us consider the same $\vartheta\in (0,1)$ that was fixed at the beginning of Section \ref{Sec Control Strategy}.
Let $\displaystyle s_\vartheta=\sum_{k=1}^\infty\,\frac{1}{k^{1+\frac{\vartheta}{2}}}$ be. Let $\tau^*_0=\tau(\frac{\rho^*_0}{2})>0$ and $M^*_0=M(\frac{\rho^*_0}{2})>0$ 
be the positive time and constant of Lemma \ref{p1}, associated to $\rho=\frac{\rho_0^*}{2}\in (0,\rho_0]$.

To obtain the proof of Theorem \ref{th2.2} it need the following Proposition \ref{P2}.

\begin{prop}\label{P2}
There exists $\ve^*_0=\ve^*_0(\frac{\rho^*_0}{2})\in (0,1)$ such that for every $\ve\in (0,\ve_0^*] 
,\,N\in \N, 
$ 
there exists
\begin{equation}\label{W^N}
W^N=\displaystyle \left\{(\tau_k)_1^{N},(w_{k})_1^{N}\right\}
\in{\cal{W^*}}(u_0) \text{ with }
 \displaystyle \tau_k 
 \leq  "\widetilde{\tau}_{k}:= \left(\frac{\ve \rho_0^*}{4M^*_0\,s_\vartheta}\right)^{\frac{2}{2+\vartheta}}\frac{1}{k}
 ,\;\;k=1,\ldots,N\; (\footnote{$\widetilde{\tau}_k$ and $\tau_k$ are defined in \eqref{tau***1} and \eqref{tau***k}, see also Remark \ref{remaining}.}). 
 \end{equation}
For every $N\in \N,\, N>n$ and $W^N
\in{\cal{W^*}}(u_0),$ as in \eqref{W^N},  we have
$\rho(W^N)\geq \frac{\rho^*_0}{2}
$ and
 \begin{equation}\label{functional J}
J^*(W^N)
\leq \sum_{l=1}^{n}\,|x^0_{l}-x^*_{l}|+c_1(\ve)\,\sum_{k=1}^{N}\frac{1}{{k}^{1+\frac{\vartheta}{2}}}
-c_2(\ve)\sum_{k=n+1}^{N}\frac{1}{k}, 
\end{equation}
where $c_1(\ve)=
\frac{\ve\rho_0^*n}{4
\,s_\vartheta}
,\; c_2(\ve)=
\left(\frac{\ve\rho_0^*}{4M^*_0
\,s_\vartheta}\right)^{\frac{2}{2+\vartheta}}.$ 
\end{prop}
Let us give the following definition.
\begin{defn}
We call {\it Separating order processing steering set of times and initial data} 
a $W^N\in {\cal{W^*}}(u_0)$ such that $\rho(W^N)\geq \frac{\rho^*_0}{2}$.
We denote with ${\cal{W^*_S}}
(u_0)$ the set of all the Separating order processing steering strategies, that is
$${\cal{W_{S}^*}}(u_0)=\left\{W^N \in {\cal{W^*}}(u_0):\;\rho(W^N)\geq \frac{\rho^*_0}{2}\right\}.$$
\end{defn}

Before proving Proposition \ref{P2}, we show how thanks this proposition, 
 we can easily obtain the proof of Theorem \ref{th2.2}.\\
{\bf Proof (of Theorem \ref{th2.2}).}
To prove Theorem \ref{th2.2} it is sufficient to show that  
\begin{equation}\label{inf}
\forall \ve>0\;\;\exists\,N_{\varepsilon}\in\N,\,W^{N_\varepsilon}\in {\cal{W_{S}^*}}(u_0)
: \qquad J^*(W^{N_\varepsilon})<\varepsilon\;(\footnote{$\displaystyle \inf_{W^N\in {\cal{W^*_S}}(u_0)}J^*(W^N)=0$.})\quad \text{ and } \quad L^{N_{\ve}}_{IS}=\{1,\ldots,n\}.
\end{equation}
This follows from Proposition \ref{P2}. Indeed, 
by contradiction,
\begin{itemize}
\item[$\star$] if $\displaystyle \inf_{W^N\in {\cal{W^*_S}}(u_0)}J^*(W^N)>0,$ 
$$\exists \bar{\ve}>0
:\,\forall N\in\N,\,\forall W^{N}\in {\cal{W_{S}^*}}(u_0) \text { we have }
 J^*(W^{N})>\bar{\ve};$$
 \item[$\star$] if $L^{N}_{IS}\neq\{1,\ldots,n\}, \forall N\in\N,$  $\quad 
  \exists j
  \not\in L^{N}_{IS},$ for every $N\in\N,$ that is
$$\exists \tilde{\ve}>0
:\,\forall N\in\N,\,\forall W^{N}\in {\cal{W_{S}^*}}(u_0) \text { we have }
 \tilde{\ve}< |\xi^{N}_{j}(T_N)-x^*_{j}|.$$
 \end{itemize}
Let us consider $\ve^*:=\min\{\bar{\ve},\tilde{\ve},\ve^*_0\},$ where $\ve_0^*\in (0,1)$ is given by Proposition \ref{P2}. In both the previous cases, 
for every $N\in\N,\,N>n,$ choose $W^N=\displaystyle \left\{(\tau_k)_1^{N},(w_{k})_1^{N}\right\}
\in{\cal{W^*}}(u_0)
$ such that
$\displaystyle
 \displaystyle \tau_k
 \leq \widetilde{\tau}_{k}= \left(\frac{\ve^* \rho_0^*}{4M^*_0\,s_\vartheta}\right)^{\frac{2}{2+\vartheta}}\frac{1}{k}
 ,\;\;k=1,\ldots,N.$ By inequality \eqref{functional J} of Proposition \ref{P2}, we obtain that, for every $N\in\N,\,N>n,$ the following inequality holds
$$ \ve^*< |\xi^N_{j}(T_{N})-x^*_{j}|\leq J^*(W^N)\leq \sum_{l=1}^{n}\,|x^0_{l}-x^*_{l}|+c_1(\ve^*)\,\sum_{k=1}^{N}\frac{1}{{k}^{1+\frac{\vartheta}{2}}}-c_2(\ve^*)\,\sum_{k=n+1}^{N}\frac{1}{k} .
$$
Keeping in mind that $\displaystyle\,\sum_{k=1}^{+\infty}\frac{1}{k}=+\infty\,,$ for $N$ enough large, from the previous inequality we have a contradiction,
then 
\eqref{inf} holds.$\qquad\qquad\qquad\qquad\qquad\qquad\qquad\qquad\qquad\qquad\qquad\qquad\qquad\diamond$
\\

Now, we give the proof of Proposition \ref{P2}. \\
{\bf Proof (of Proposition \ref{P2}).}
Let $(\s_k)_1^N\in \R^N_+$ be. Given $
N\in\N$ and a generic $W^N=\displaystyle \left\{(\tau_k)_1^{N},(w_{k})_1^{N}\right\}
\in{\cal{W^*}}(u_0),$ 
let us recall that the curves of sign change 
associated to 
 $W^N$ (see Section \ref{OPS}) are defined in the following way
$$\xi_l^{W}(t)=
\xi^{k
}_l
(t
),\qquad \forall t\in [S_{k},T_k], 
\;\; k=1,\ldots,N, \;l=0,\ldots,n+1,
$$
where all the {curves}
$\xi^{k
}_l,  \;l=0,\ldots,n+1,
$ 
are defined on $[S_k,T_k]$ and are associated to
the initial state $w_{k}.$ 
Let us recall that in Section \ref{OPS}, for every $k=2,\ldots,N,$ we have defined
$x^{k-1}_j:=\xi^{k-1}_j(T_{k-1})
=\xi^{k}_j(S_{k}), \;j=0,\ldots,n+1$
 and $\displaystyle \rho_{k-1}=\min_{l=0,\ldots,n}\{x^{k-1}_{l+1}-x^{k-1}_l\}.$
We recall that by Lemma \ref{p1} (see also the properties of $\big\{\tau_k, w_k\big\},\,k=1,\ldots,N,$ in Section \ref{OPS}) 
\begin{align}\label{Mk}
\forall\rho\in(0,\rho_{k-1}]\;\;
 \exists\,\overline{\tau}_k=\overline{\tau}_k(\rho)>0,\; &M_k=M_k(\rho)>0,\,
\xi_l^{k}\in C^{1+\frac{\vartheta}{2}}([S_k,T_k]),\;l=0,...,n+1,\;\nonumber\\
\\
\text{  such that  } \;\;\;\;\;\;&
\|\xi^{k}_l
\|_{C^{1+\frac{\vartheta}{2}}([S_k,T_k])}\leq M_{k}.\nonumber  
\end{align}
By the properties of $\displaystyle \left\{\tau_k, w_{k}\right\},$ keeping in mind  
Remark \ref{rem steer},
we observe that
\begin{equation}\label{steer}
\dot{\xi^k_l}
(S_k) \; = 
\mu_k(x^*_l-x_j^0), \;\;\;\;l=0,\ldots,n+1.
\end{equation}
{\bf Step 1:} {\it Some preliminary evaluations.}\\
 {\it Gap estimate.
} 
Let $k=1,\ldots,N$ and $l=0,
\ldots,n
$ be.
We note that 
\begin{align}\label{est gap}
\hspace{-0.25cm}
\xi^k_{l+1}(t)-\xi^k_l(t)=&\xi^k_{l+1}(t)-x^{k-1}_{l+1}+x^{k-1}_{l+1}-x^{k-1}_{l}+x^{k-1}_{l}-\xi^{k}_l(t)\nonumber\\
=&  x^{k-1}_{l+1}-x^{k-1}_{l}+\int_{S_k}^t\left(\dot{\xi}_{l+1}^k(s)-\dot{\xi}_{l}^k(s)\right)\,ds\nonumber\\
=&  x^{k-1}_{l+1}-x^{k-1}_{l}
+\int_{S_k}^t\left(\dot{\xi}_{l+1}^k(s)-\dot{\xi}_{l+1}^k(S_k)\right)\,ds
-\int_{S_k}^t\left(\dot{\xi}_{l}^k(s)-\dot{\xi}_{l}^k(S_k)\right)\,ds
\nonumber\\
&\qquad\qquad\qquad\qquad\qquad
+\left(\dot{\xi}_{l+1}^k(S_k)-\dot{\xi}_{l}^k(S_k)\right)\left(t-S_k\right),\,\;\forall t\in [S_k,T_k]
\,.
\end{align}
By \eqref{Mk}, for every $j=0,\ldots,n+1,$ 
we deduce that
\begin{equation}\label{hold1}
\int_{S_k}^t\left|\dot{\xi^k_j}(s)-\dot{\xi^k_j}(S_k)\right|\,ds\leq M_k\int_{S_k}^t\,\left(s-S_k\right)^{\frac{\vartheta}{2}}\,ds
\leq M_k\left(t-S_k\right)^{1+\frac{\vartheta}{2}}, \;\;\quad\forall t\in [S_k,T_{k}].
\end{equation}
By \eqref{hold1} 
we have
\begin{equation}\label{gap0}
\! \left|\int_{S_k}^t\left(\dot{\xi}_{l+1}^k(s)-\dot{\xi}_{l+1}^k(S_k)\right)\,ds-\int_{S_k}^t\left(\dot{\xi}_{l}^k(s)-\dot{\xi}_{l}^k(S_k)\right)\,ds\right|\leq 
2M_k\left(t-S_k\right)^{1+\frac{\vartheta}{2}}\!\!\!\!,
\forall t\in [S_k,T_k].
\end{equation}
By \eqref{steer}, \eqref{est gap} and \eqref{gap0}, 
for every $ t\in [S_k,T_k],$ we deduce
\begin{equation}\label{for step4}
\xi^k_{l+1}(t)-\xi^k_l(t)\geq x^{k-1}_{l+1}-x^{k-1}_{l}-2M_{k}
\tau_k^{1+\frac{\vartheta}{2}} +\left(\mu_k(x_{l+1}^*-x_{l+1}^0)-\mu_k(x_{l}^*-x_{l}^0)
\right)\left(t-S_k\right).
\end{equation}
{\it Distance evaluation.}
Let $k=1,\ldots,N$ and $l=1,\ldots,n$ be.
By  
\eqref{steer} we have
\begin{align}\label{target dist}
\xi^k_l(t)-x^*_l=&\xi^k_l(t)-x^{k-1}_l+x^{k-1}_l-x^*_l=x^{k-1}_l-x^*_l+\int_{S_k}^t\dot{\xi^k_l}(s)ds\nonumber\\
=&x^{k-1}_l-x^*_l+\dot{\xi^k_l}(S_k)\left(t-S_k\right)+\int_{S_k}^t\left(\dot{\xi^k_l}(s)-\dot{\xi^k_l}(S_k)\right)ds\nonumber\\
=& x^{k-1}_l-x^*_l+\mu_k(x^*_l-x_l^0)\,\left(t-S_k\right)+\int_{S_k}^t\left(\dot{\xi^k_l}(s)-\dot{\xi^k_l}(S_k)\right)ds,\,\forall t\in[S_{k},T_k].
\end{align}
{\bf Step 2:}{ \it Uniform time $\tau_0^*$ and constant $M^*_0$.} Let $k=1,\ldots,N$.
Without loss of generality, we can suppose, by an induction argument, that we have already proved 
that
\begin{equation}\label{ind assum}
 \displaystyle\rho_{h-1}=\min_{l=0,\ldots,n}\{x^{h-1}_{l+1}-x^{h-1}_l\}\geq\frac{\rho^*_0}{2}, \text{ for every } h=1,\ldots,k. \;(\footnote{Let us recall that by definition $ \displaystyle\rho_{0}=\min_{l=0,\ldots,n}\{x^{0}_{l+1}-x^{0}_l\}\geq\frac{\rho^*_0}{2}$.})
 \end{equation}
  Then, choosing $\rho=\frac{\rho^*_0}{2}\in(0,\rho_{h-1}],$ $h=1,\ldots,k$ in \eqref{Mk}, we obtain 
  \begin{equation}\label{Mh}
  M_h=M_h\Big(\frac{\rho_0^*}{2}\Big)=M^*_0,\;\;\;\overline{\tau}_h=\overline{\tau}_h\Big(\frac{\rho_0^*}{2}\Big)=\tau_0^*,\;\forall h=1,\ldots,k.
    \end{equation}
Define $\ve^*_0=\min\left\{\frac{4M^*_0s_{\vt}}{\rho_0^*}{\tau_0^*}^{1+\frac{\vt}{2}},1\right\}$, fix $\ve\in(0,\ve_0^*],$ and let 
 \begin{equation}\label{tau uniform0}
 \displaystyle \tau_h 
 \leq  \widetilde{\tau}_{h}= \widetilde{\tau}_{h}(\ve):= \left(\frac{\ve \rho_0^*}{4M^*_0\,s_\vartheta}\right)^{\frac{2}{2+\vartheta}}\frac{1}{h}
 ,\;\;h=1,\ldots,k\,.\; (\footnote{$\widetilde{\tau}_h$ and $\tau_h$ are defined in \eqref{tau***1} and \eqref{tau***k}, see also Remark \ref{remaining}.})
 \end{equation}
 So, we have
  \begin{equation}\label{tau uniform}
  \widetilde{\tau}_{h}(\ve)\leq \tau_0^*,\;\;h=1,\ldots,k. 
   \end{equation}
In the final Step.5, by \eqref{for step4}, \eqref{Mh}-\eqref{tau uniform}, and a technical proof, for every $l=0,\ldots,n$ and $k=1,\ldots,N,$ 
we will prove the following alternative inequalities 
\begin{align}\label{gapk}
\xi^{k}_{l+1}(t)-\xi^{k}_l(t)
\geq& x^{k-1}_{l+1}-x^{k-1}_{l}
-2
M_0^* \tau_k^{1+\frac{\vartheta}{2}}
,\;\;\,\:\qquad\forall t\in 
[S_k,T_k],\nonumber\\ 
\text{ or }&\\
\xi^{k}_{l+1}(t)-\xi^{k}_l(t)\geq &\rho^*_0-\ve \frac{\rho^*_0}{4}\geq \frac{3}{4}\rho^*_0,  
\qquad\qquad\qquad\quad\;\forall t\in [S_k,T_k]
.\nonumber
\end{align}
{\bf Step 3:}{ \it Let us prove $\rho(W^N)\geq\frac{\rho^*_0}{2}.$} In this step, we will prove that $\xi^{k}_{l+1}(t)-\xi^{k}_l(t)\geq \frac{\rho^*_0}{2},\,\forall t\in 
[S_k,T_k],\:k=1,\ldots,N.$
By \eqref{tau uniform0}, \eqref{gapk} for k=1, we have
\begin{align*}
\min_{l=0,\ldots,n}\min_{t\in [S_1,T_{1}]}\left(\xi^{1}_{l+1}(t)-\xi^{1}_l(t)\right)\geq &
\min\left\{\min_{l=0,\ldots,n}\left(x^{0}_{l+1}-x^{0}_{l}\right)
-2M^*_0 \tau_1^{1+\frac{\vartheta}{2}}, \frac{3}{4}\rho^*_0\right\}\nonumber\\
\geq &
\min\left\{\rho_0^*
-2M^*_0\,\frac{\ve\rho_0^*}{4M^*_0\,s_\vartheta}, \frac{3}{4}\rho^*_0\right\} 
\nonumber\\= &\min\left\{\rho_0^*\left(1-\frac{\ve}{2s_\vartheta}\right), \frac{3}{4}\rho^*_0\right\}\geq\frac{\rho_0^*}{2}.
\end{align*}
By \eqref{tau uniform0}, \eqref{gapk} for $k=2$ and $k=1,$ 
we obtain
\begin{align*}
\hspace{-0.5cm}\min_{t\in [S_1,T_{1}]
}\min_{l=0,\ldots,n}\left(\xi^{2}_{l+1}(t)-\xi^{2}_l(t)\right)\geq &
\min\left\{\min_{l=0,\ldots,n}\left(x^{1}_{l+1}-x^{1}_{l}\right)
-2M^*_0 \tau_2^{1+\frac{\vartheta}{2}}, \frac{3}{4}\rho^*_0\right\}\nonumber\\
\geq &\min\left\{\min_{l=0,\ldots,n}\left(x^{0}_{l+1}-x^{0}_{l}
\right)
-2M^*_0 \left(\tau_1^{1+\frac{\vartheta}{2}}+\tau_2^{1+\frac{\vartheta}{2}}\right), \frac{3}{4}\rho^*_0\right\}\nonumber\\
\geq&\min\left\{\rho_0^*-2M^*_0\,\frac{\ve\rho_0^*}{4M^*_0\,s_\vartheta} \left(1+\frac{1}{2^{1+\frac{\vartheta}{2}}}\right), \frac{3}{4}\rho^*_0\right\}\nonumber\\
=&\min\left\{\rho_0^*\left(1-\frac{1+\frac{1}{2^{1+\frac{\vartheta}{2}}}}{2s_\vartheta}\ve\right), \frac{3}{4}\rho^*_0\right\}\geq\frac{\rho_0^*}{2}.
\end{align*}
In general, in the case $N\geq3,$ 
 for every $k=3,\ldots, N,$ by \eqref{tau uniform0} and \eqref{gapk} we deduce
\begin{align*}
\min_{l=0,\ldots,n}&\min_{t\in [S_k,T_k]}\left(\xi^{k}_{l+1}(t)-\xi^{k}_l(t)\right)\geq  
\min\left\{\min_{l=0,\ldots,n}\left(x^{k-1}_{l+1}-x^{k-1}_{l}
\right)
-2M^*_0 \tau_k^{1+\frac{\vartheta}{2}}, \frac{3}{4}\rho^*_0\right\}\nonumber\\
\geq & 
\min\left\{\min_{l=0,\ldots,n}\left(x^{0}_{l+1}-x^{0}_{l}
\right)
-2M^*_0 \left(\tau_1^{1+\frac{\vartheta}{2}}+\tau_2^{1+\frac{\vartheta}{2}}+\ldots+\tau_k^{1+\frac{\vartheta}{2}}\right), \frac{3}{4}\rho^*_0\right\}\nonumber\\
\geq &\min\left\{\rho_0^*-2M^*_0\,\frac{\ve\rho_0^*}{4M^*_0\,s_\vartheta} \left(1+\frac{1}{2^{1+\frac{\vartheta}{2}}}+\ldots+\frac{1}{{k}^{1+\frac{\vartheta}{2}}}\right), \frac{3}{4}\rho^*_0\right\}\nonumber\\
=&\min\left\{\rho_0^*\left(1-\frac{1+\frac{1}{2^{1+\frac{\vartheta}{2}}}+\ldots+\frac{1}{{k}^{1+\frac{\vartheta}{2}}}}{2s_\vartheta}\ve\right), \frac{3}{4}\rho^*_0\right\}\geq\frac{\rho_0^*}{2}.
\end{align*}
Thus, it follows
$$\rho(W^N):=\min_{l=0,\ldots,n}\min_{t\in 
\mE
}
\left\{\xi_{l+1}^{W}(t)-\xi_l^{W}(t)\right\}=\min_{l=0,\ldots,n}\min_{k=1,\ldots,N}\min_{t\in [S_k,T_k]}
\left(\xi^{k}_{l+1}(t)-\xi^{k}_l(t)\right)\geq\frac{\rho_0^*}{2}.$$
{\bf Step 4:} {\it Proof of \eqref{functional J}.} For $N\in\N,\,N>n,$ let $W^N=
\displaystyle \left\{(\tau_k)_1^{N},(w_{k})_1^{N}\right\}\in{\cal{W^*}}(u_0)$, as in \eqref{W^N}. 
Let $k=1,\ldots,N$ and $l=1,\ldots,n$. 
Since in Step. 3 we have proved that $\rho(W^N)\geq\frac{\rho^*_0}{2},$
we note that inequality \eqref{hold1} holds with $M_k=M_k(\frac{\rho^*_0}{2})=M_0^*.$ 
Then, 
\begin{itemize}
\item[$\star$] 
if
$l\in
L^{k-1}_{IS},$ keeping in mind that $\mu_k(x^*_l-x_l^0)=0,$ by \eqref{target dist} and \eqref{hold1} with $M_k=M_0^*,$ we obtain
\begin{equation}\label{est0}
|\xi^k_{l}(t)-x^*_{l}|\leq |x^{k-1}_{l}-x^*_{l}|+M^*_0\,\left(t-S_k\right)^{1+\frac{\vartheta}{2}},
\qquad \qquad \qquad \;\;\quad\forall t\in[S_{k},T_k]\,;
\end{equation}
\item[$\star$] if $l\in\{1,\ldots,n\}\backslash L^{k-1}_{IS},$ keeping in mind that $x^{k-1}_l\neq x_l^*,$ we have
\begin{align}\label{est1}
x^{k-1}_{l}-x^*_{l}+\mu_k(x^*_{l}-x_{l}^0)\,\left(t-S_k\right)=& (x^{k-1}_{l}-x^*_{l})\left(1+\frac{\mu_k(x^*_{l}-x_{l}^0)}{x^{k-1}_{l}-x^*_{l}}\left(t-S_k\right)\right)\nonumber\\
=&(x^{k-1}_{l}-x^*_{l})\left(1-\frac{t-S_k}{|x^{k-1}_{l}-x^*_{l}|}\right),\; \forall t\in[S_{k},T_k].
\end{align}
Thus, by \eqref{target dist}, \eqref{est1} and \eqref{hold1} with $M_k=M_0^*,$ keeping in mind that by the definition of $\tau_k$ (see \eqref{tau***1} and \eqref{tau***k}) 
it follows that $\displaystyle \tau_k\leq \min_{l\in\{1,\ldots,n\}\backslash L^{k-1}_{IS}}|x^{k-1}_{l}-x^*_{l}|,$ 
we deduce
\begin{align}\label{tau stop}
|\xi^k_{l}(t)-x^*_{l}|&\leq |x^{k-1}_{l}-x^*_{l}|\left|1-\frac{t-S_k}{|x^{k-1}_{l}-x^*_{l}|}\right|+M^*_0\,\left(t-S_k\right)^{1+\frac{\vartheta}{2}}
\nonumber\\
&=|x^{k-1}_{l}-x^*_{l}|-\left(t-S_k\right)+M^*_0\,\left(t-S_k\right)^{1+\frac{\vartheta}{2}},
\; \qquad \qquad \;\;\forall t\in
\Big[S_k,T_k\Big].
\end{align}
\end{itemize}
Let us set $\displaystyle J_k(t):=\sum_{l=1}^{n}\,|\xi^k_{l}(t)-x^*_{l}|, \,t\in [S_{k},T_k].$\\
By \eqref{est0} and \eqref{tau stop}, for every $k=1,\ldots,N,$ for every $t\in \left[S_k,T_k
\right],$ 
we obtain
\begin{align}\label{57}
J_k(t)=&
\sum_{l\in\{1,\ldots,n\}\backslash L^{k-1}_{IS}}\,|\xi^k_{l}(t)-x^*_{l}|+\sum_{
l\in L^{k-1}_{IS}
}\,|\xi^k_{l}(t)-x^*_{l}|\nonumber\\ 
\leq &\sum_{l\in\{1,\ldots,n\}\backslash L^{k-1}_{IS}}\,|x^{k-1}_{l}-x^*_{l}|-\left(n-card(L^{k-1}_{IS})\right)\,\left(t-S_k\right)\nonumber\\
+&\left(n-card(L^{k-1}_{IS})\right)\,M^*_0\,\left(t-S_k\right)^{1+\frac{\vartheta}{2}}
+\sum_{l\in L^{k-1}_{IS}}|x^{k-1}_{l}-x^*_{l}|+card(L^{k-1}_{IS})\,M^*_0\,\left(t-S_k\right)^{1+\frac{\vartheta}{2}}
\nonumber\\
=&J_{k}(S_k)-\left(n-card(L^{k-1}_{IS})\right)\,\left(t-S_k\right)+n\,M^*_0\,\left(t-S_k\right)^{1+\frac{\vartheta}{2}}
\nonumber\\
\leq &  J_{k}(S_k)+n\,M^*_0\,\,\tau_k^{1+\frac{\vartheta}{2}}-\left(t-S_k\right),\qquad\qquad\qquad\qquad\qquad\qquad \forall t\in \left[S_k,T_k
\right]
.
\end{align}
Now, keeping in mind that $J_h(S_h)=J_{h-1}(T_{h-1}),\, h=2,\ldots,N,$ by \eqref{57} we can deduce that
\begin{align}\label{J_k}
J_1(t)=&\sum_{l=1}^{n}\,|\xi^1_{l}(t)-x^*_{l}|\leq J_1(S_1)+n\,M^*_0\,\,\tau_{1}^{1+\frac{\vartheta}{2}}-\left(t-S_1\right) 
,\;\forall t\in [S_1,T_{1}],\nonumber\\
J_2(t)=&\sum_{l=1}^{n}\,|\xi^2_{l}(t)-x^*_{l}|\leq J_2(S_2)+n\,M^*_0\tau_{2}^{1+\frac{\vartheta}{2}}-\left(t-S_2\right) 
= J_1(T_1)+n\,M^*_0{\tau_{2}}^{1+\frac{\vartheta}{2}}\!\!\!-\left(t-S_2\right)
\nonumber\\
\leq & J_1(S_1)+n\,M^*_0\left(\tau_{1}^{1+\frac{\vartheta}{2}}+\tau_{2}^{1+\frac{\vartheta}{2}}\right)-\left(\tau_1+t-S_2\right),\;\;\;\quad\forall t\in [S_2,T_{2}],\nonumber\\
J_k(t)=&\sum_{l=1}^{n}\,|\xi^k_{l}(t)-x^*_{l}|\leq J_k(S_k)+n\,M^*_0\,\,\tau_k^{1+\frac{\vartheta}{2}}-\left(t-S_k\right)
\nonumber
\\\leq&
J_1(S_1)+n\,M^*_0\,\left(\tau_{1}^{1+\frac{\vartheta}{2}}+\ldots+\tau_k^{1+\frac{\vartheta}{2}}\right)-
\left(\tau_1+\ldots+\tau_{k-1}+t-S_k\right),\nonumber\\ 
&\qquad\qquad\qquad\qquad\qquad\qquad\qquad\qquad\qquad\qquad\qquad\;\,\forall   t\in [S_k,T_k],\,\forall k=3,\ldots,N.
\end{align} 
By \eqref{J_k}, keeping in mind that $\displaystyle J_1(S_1)= \sum_{l=1}^{n}\,|x^0_{l}-x^*_{l}|$ and $L^{N-1}_{IS}\subseteq\{1,\ldots,n\},$
 by the definition of $\widetilde{\tau}_k$ and $\tau_k$ (see \eqref{tau***1} and \eqref{tau***k}) and Remark \ref{remaining}, inequality \eqref{functional J} follows. Indeed, 
 \begin{align*}
 J^*(W^N)=&\sum_{l=1}^{n}\,|\xi^W_{l}(T_{N})-x^*_{l}|=J_N(T_{N})\leq \sum_{l=1}^{n}\,|x^0_{l}-x^*_{l}|+nM^*_0\sum_{k=1}^{N}\tau_k^{1+\frac{\vartheta}{2}}-\,\sum_{k=1}^{N}\tau_k\\\leq& \sum_{l=1}^{n}\,|x^0_{l}-x^*_{l}|+\frac{\ve{\rho}^*_0n}{4s_\vt}\sum_{k=1}^{N}\frac{1}{{k}^{1+\frac{\vartheta}{2}}}
-\,\sum_{\{k\,|\,\Theta_k\not=\emptyset\}
}
\tau_{k}-\left(\frac{\ve\rho_0^*}{4M^*_0\,s_\vartheta}\right)^{\frac{2}{2+\vartheta}}\!\!\!\!\!\!\!\!\!\sum_{\{k\,|\,\Theta_k=\emptyset\}
}\frac{1}{k}\\
\leq & \sum_{l=1}^{n}\,|x^0_{l}-x^*_{l}|+\frac{\ve{\rho}^*_0n}{4s_\vt}\sum_{k=1}^{N}\frac{1}{{k}^{1+\frac{\vartheta}{2}}}
-\left(\frac{\ve\rho_0^*}{4M^*_0\,s_\vartheta}\right)^{\frac{2}{2+\vartheta}}\!\!\!\!\!\sum_{k=n+1}^N\frac{1}{k}.
\end{align*}
{\bf Step 5:} {\it Proof of \eqref{gapk}.}  In this final step, assuming that 
\eqref{ind assum} holds, we prove that inequality \eqref{for step4} implies the two alternative inequalities in \eqref{gapk}. 
Let us set
$$A^k_l:=
\mu_k(x_{l+1}^*-x_{l+1}^0)-\mu_k(x_{l}^*-x_{l}^0),\;\;\, \forall k=1,\ldots,N,\;\forall l=0,\ldots,n.$$
Let $k=1,\ldots,N$ and $l=0,\ldots,n$.\\
 We note that if $l,l+1\in L^{k-1}_{IS}$ it follows $A^k_l=0,$ thus we obtain the first of the two alternatives in \eqref{gapk}.
  Instead, if at least one between $l$ and $l+1$ is not an inactive index (that is $l,l+1\not\in L^{k-1}_{IS}$ or $l\in L^{k-1}_{IS},l+1\not\in L^{k-1}_{IS}$ or $l\not\in L^{k-1}_{IS},l+1\in L^{k-1}_{IS}$), since the conditions $$x^{k-1}_l<x^{k-1}_{l+1}\;\;\;\text{ and } \;\;\;x^*_l<x^*_{l+1}, \;\;\; l=0,\ldots,n,$$
hold, we have to distinguish six possible configurations to calculate $A^k_l$, that are presented in the following table:\\
\begin{table}[ht]
\caption{}
\centering
\begin{tabular}{|c|c|c|}
\hline
Configurations of the points &  i) &ii)\\
\hline
 a & $x_l^{k-1}<x_{l+1}^{k-1}\leq x_{l}^*<x_{l+1}^*$ & $x_l^*<x_{l+1}^*\leq x_{l}^{k-1}<x_{l+1}^{k-1}$\\
\hline
 b & $x_l^{k-1}\leq x_{l}^*\leq x_{l+1}^{k-1}\leq x_{l+1}^*$ & $x_l^*\leq x_{l}^{k-1}\leq x_{l+1}^*\leq x_{l+1}^{k-1}$ \\
\hline
 c & $x_l^{k-1}\leq x_{l}^*<x_{l+1}^*\leq x_{l+1}^{k-1}$ & $x_l^*\leq x_{l}^{k-1}<x_{l+1}^{k-1}\leq x_{l+1}^*$ \\
\hline
\end{tabular}\\
\label{table 1}
\end{table}

\noindent Thus, in the following table we can compute $A^k_l$ in any configuration 
in the 3 cases: \\
$l,l+1\not\in L^{k-1}_{IS};$\; $l\in L^{k-1}_{IS}$ and $l+1\not\in L^{k-1}_{IS};$ \;$l\not\in L^{k-1}_{IS}$ and $l+1\in L^{k-1}_{IS}.$  
\begin{table}[ht]
\centering
\caption{\quad\qquad\quad\;\;\;\quad$l,l+1\not\in L^{k-1}_{IS}$$\quad\quad l\in L^{k-1}_{IS},l+1\not\in L^{k-1}_{IS}$\quad$l\not\in L^{k-1}_{IS},l+1\in L^{k-1}_{IS}$}
\hspace{2.5cm}
\begin{tabular}{|c|c|c|}
\hline
$A^k_l$ 
&  i) &ii)\\
\hline
 a &0 &0\\
\hline
 b &0 &0 \\
\hline
 c &-2 $^\clubsuit$&2 \\
\hline
\end{tabular}
\hspace{0.5cm}
\begin{tabular}{|c|c|c|}
\hline
$A^k_l$ 
&  i) &ii)\\
\hline
 a &1 &-1$^\diamondsuit$\\
\hline
 b &1 & -1$^\diamondsuit$\\
\hline
 c &-1$^\diamondsuit$ &1 \\
\hline
\end{tabular}
\hspace{0.5cm}
\begin{tabular}{|c|c|c|}
\hline
$A^k_l$ 
&  i) &ii)\\
\hline
 a &-1$^\spadesuit$ &1\\
\hline
 b &-1$^\spadesuit$ &1 \\
\hline
 c &-1$^\spadesuit$ &1 \\
\hline
\end{tabular}
\label{table 2}
\end{table}

\noindent We can observe that when $A^k_l\geq 0$ (in Table \ref{table 2}), that is in the cases and in the configurations without 
superscript symbols (without $\clubsuit$ or $\diamond$
or $\spadesuit$), by \eqref{for step4} we easily obtain the first of the two alternative inequalities in \eqref{gapk}.\\

Before analyzing the cases with superscript symbols,
we remark that 
if
$j\in
L^{k-1}_{IS},$
there exists $h_j,\,0\leq h_j\leq k-1$ such that $x^{h_j}_j=x^*_j,$ then
 keeping in mind that $\mu_{h_j}(x^*_j-x_j^0)
 =\ldots=\mu_k(x^*_j-x_j^0)=0,$ thus, by \eqref{target dist}, \eqref{hold1}, \eqref{ind assum}-\eqref{tau uniform}, we obtain 
\begin{align}\label{contr}
|\xi^k_{j}(t)-x^*_{j}|\leq & |x^{k-1}_{j}-x^*_{j}|+M^*_0\,\tau_{k}^{1+\frac{\vartheta}{2}}\leq |x^{k-2}_{j}-x^*_{j}|+M^*_0\,\left(\tau_{k-1}^{1+\frac{\vartheta}{2}}+\tau_{k}^{1+\frac{\vartheta}{2}}\right)\nonumber\\
\leq &|x^{h_j}_{j}-x^*_{j}|+M^*_0\,\left(\tau_{h_j+1}^{1+\frac{\vartheta}{2}}+\ldots+\tau_{k-1}^{1+\frac{\vartheta}{2}}+\tau_{k}^{1+\frac{\vartheta}{2}}\right)\nonumber\\
=&M^*_0\sum_{s=h_j+1}^k \tau_{s}^{1+\frac{\vartheta}{2}} \leq M^*_0\frac{\ve \rho_0^*}{4M^*_0\,s_\vartheta}\sum_{s=h_j+1}^k \frac{1}{k^{1+\frac{\vartheta}{2}}}\leq \ve\frac{\rho_0^*}{4},\qquad \qquad\forall t\in[S_{k},T_k],\,
\end{align}
so, 
the point of sign change with index $j\in L^{k-1}_{IS}$
 remain forever near the corresponding target point of sign change already reached.\\
 
Now, we can analyze the following 3 cases:
\begin{description}
\item[$\clubsuit$] In the case of $l,l+1\not\in L^{k-1}_{IS},$ by the configuration $i.c),$ for every $t\in[S_k,T_{k}]$ we have
$$\!\!\!\!\!\!\!\!\!\!\!\!\! \!\!\!\!  x_l^{k-1}\leq\xi_l^k(t)\leq x_{l}^*<x_{l+1}^*\leq\xi_{l+1}^k(t)\leq x_{l+1}^{k-1}\Longrightarrow\xi_{l+1}^k(t)-\xi_l^k(t)\geq x_{l+1}^*-x_{l}^*\geq \rho^*_0
,$$
from which the second option of \eqref{gapk} follows.
\item[$\diamondsuit$]
In the case of $l\in L^{k-1}_{IS},l+1\not\in L^{k-1}_{IS},$ we analyze the following 3 configurations: 
\begin{description}
\item[ii.a)] By \eqref{contr}, we deduce $x_l^*<x_{l+1}^*\leq x_{l}^{k-1}\leq x_{l}^{*}+\frac{\ve\rho}{4},$ then $\rho^*_0\leq x_{l+1}^*-x_l^*\leq x_{l}^{k-1}-x_l^*\leq \frac{\ve\rho^*_0}{4}\leq \frac{\rho^*_0}{4},$ from which a contradiction follows, so this configuration 
is not admissible.
\item[ii.b)] By \eqref{contr}, we deduce 
\begin{align*}
 \xi^k_{l}(t)\leq x_l^*+ \frac{\ve\rho^*_0}{4}
\leq x_{l+1}^*\leq\xi_{l+1}^k(t)\leq x_{l+1}^{k-1} \Longrightarrow\xi_{l+1}^k(t)-\xi_l^k(t)\geq & x_{l+1}^*-\left(x_{l}^*+\frac{\ve\rho^*_0}{4}\right)\\
\geq& \rho^*_0-\frac{\ve\rho^*_0}{4}\geq \frac{3}{4}\rho^*_0,
\end{align*}
from which the second option of \eqref{gapk} follows.
\item[i.c)] It is similar to the configuration $ii.b)$. 
\end{description}
\item[$\spadesuit$]
In the case of $l\not\in L^{k-1}_{IS},l+1\in L^{k-1}_{IS},$ we analyze the following 3 configurations:
\begin{description}
\item[i.a)] It is similar to $ii.a)$ of the case $\diamondsuit$.
\item[i.b)] By \eqref{contr}, we deduce 
\begin{align*}
x_l^{k-1}\leq \xi^k_{l}(t)\leq x_l^*
\leq x_{l+1}^*-\frac{\ve\rho^*_0}{4}%
\leq \xi_{l+1}^k(t)
\Longrightarrow\xi_{l+1}^k(t)-\xi_l^k(t)\geq & \left(x_{l+1}^*-\frac{\ve\rho^*_0}{4}\right)-x_{l}^*\\
\geq& \rho^*_0-\frac{\ve\rho^*_0}{4}\geq \frac{3}{4}\rho^*_0,
\end{align*}
from which the second option of \eqref{gapk} follows.
\item[i.c)] It is similar to the configuration $i.b)$. \qquad\qquad\qquad\qquad\qquad\qquad\qquad\qquad\qquad$\diamond$
\end{description}
\end{description}

%

\section{
Proof of Theorem \ref{th preserving}}\label{regular section} 
In 
this section, we prove Theorem \ref{th preserving}. Without loss of generality, we can reformulate problem \eqref{u pres} with a generic time interval $(0, T)$, in the following way
\begin{equation}\label{u pres proof}
   \begin{cases}
\quad   u_t \; = \; u_{xx} \; + \; v(x,t)  u \; + \; f(u)
&\quad  {\rm in} \;\;\;Q_{T}= (0, 1) \times (0, T),  
\\
\quad u (0,t) = u (1,t) = 0,
&\quad t \in (0, T),
\\
\quad u\:\mid_{t = 0} \; = u_{in} + r_{in},
\end{cases}
\end{equation}
where $ u_{in},\,r_{in}\in 
H^1_0(0,1),$ and $u_{in}$ have exactly $n$ points of sign change at $ x_{l} \in (0,1),\, 
$ 
with $0=:x_0<x_l<x_{l+1}\leq x_{n+1}:=1,\;\; l=1,\ldots,n$.\\
Throughout this section, we represent the solution 
of \eqref{u pres proof} as the sum of two functions $ w (x,t) $ and $h (x,t)$, which solve the following problems in $Q_T$
\begin{equation}\label{w}
   \begin{cases}
\quad   w_t \; = \; w_{xx} \; + \; v(x,t)  w \; + \; f(w)
\\
\quad w (0,t) = w (1,t) = 0
\\
\quad w \:\mid_{t = 0} \; = u_{in}, 
\end{cases}
\;
   \begin{cases}
\quad   h_{t} \; = \; h_{xx} \; + \; v(x,t)  h \; + \; (f(w+h) - f(w))
\\
\quad h (0,t) = h (1,t) = 0
\\
\quad h \:\mid_{t = 0} \; = r_{in}. 
\end{cases}\!\!\!\!\!
\end{equation}
In this section, we denote the target state by $ \overline{u} \in H^1_0(0,1)$ instead of the specific
$ w_k$ introduced in the statement of Theorem \ref{th preserving}.

\begin{lem}\label{lem preserving} 
Let $\overline{u}\in 
H^1_0(0,1)$ have the same points of sign change of $u_{in},$ in the same order of sign change. Let us suppose that
\begin{equation}\label{ass}
\exists \nu>0:\;\nu\leq\frac{\overline{u}(x)}{u_{in}(x)}< 1, \qquad  
\forall \,x\in \left(0,1\right)\backslash \bigcup_{l=1}^n\left\{x_l\right\}. 
\end{equation}
Then, for every $\eta>0$ there exist a sufficiently small time $T=T(\eta,u_{in}, \overline{u})>0$ 
and a piecewise static bilinear control $v=v(\eta,u_{in}, \overline{u})
\in L^\infty(Q_T
)$ such that
\begin{equation}\label{ineq pres}
\|u (\cdot, T) - \overline{u}(\cdot)\|_{L^2(0,1)}\leq \eta+\sqrt{2}\|r_{in}\|_{L^2(0,1)},
\end{equation}
 where $u$ is the solution of \eqref{u pres proof} on $Q_T.$
\end{lem}

\n
{\bf Proof.}
Here, we adapt the proof introduced in Section 2 of \cite{CanKh} from the linear case to the semilinear problem \eqref{u pres proof}. We consider
the following function defined on $[0,1]$
\begin{equation*}
v_0 (x) \; = \;  \left\{ \begin{array}{ll}
\ln \left( \frac{\overline{u} (x)}{u_{in} (x) } \right), \;\;\;\; &  {\rm for } \;\; x \neq  0, 1, x_l,  \; l = 1\ldots, n   \\ 
0, \;\;\;\; &  {\rm for} \;\; x =  0, 1, x_l,  \; l = 1\ldots, n \\
\end{array}
\right.  
\end{equation*}

\bigskip
\n
By \eqref{ass} we note that 
$v_0\in L^\infty (0,1)$  
and
$
v_0 (x) \; \leq 0,
\text{ for every } x\in [0,1].$ 
We select the bilinear control
$$
v (x,t) \; := \; \frac{1}{T} v_0 (x).
$$
Then, let us represent the solution $ u $ of \eqref{u pres proof}, associated to the previous choice of the coefficient $v,$ as a sum of two functions $ w (x,t) $ and $h (x,t)$, which solve the two problems introduced in \eqref{w}, respectively.\\
{\bf Step 1:} {\it Representation formula for $ w (\cdot, T) $.}
For every fixed $\bar{x}\in (0,1),$ let us consider the non-homogeneous first-order ODE
 $ w^\prime(\bar{x},t)=\frac{v_0 (\bar{x})}{T}  w(\bar{x},t) \; +\big( w_{xx}(\bar{x},t)+f(w(\bar{x},t))\big),\; t\in (0,T),$ associated to the first problem in \eqref{w}.
Then, we easy deduce that the corresponding solution $w$ to the first problem in \eqref{w} admits the following representation
$$
w (x,t) \; = \; e^{v_0 (x) \frac{t}{T}} u_{in} (x)  \; + \; 
\int_0^t e^{v_0 (x) \frac{(t - \tau)}{T}} (w_{xx} (x, \tau) + f (w(x, \tau)))d \tau,\,\;\;\;\forall (x,t)\in Q_T,
$$
and, at time $ t = T$ we have
\begin{equation}\label{WT}
w (x,T) \; = \; \overline{u} (x)  \; + \int_0^T e^{v_0 (x) \frac{(T - t)}{T}} (w_{xx} (x, \tau) + f (w(x,\tau))) dt,\quad \forall x\in(0,1).
\end{equation}
Let us show that the integral in the right-hand side of (\ref{WT}) tends to zero in $ L^2 (0,1) $ as $ T \rightarrow 0+$, which would mean that $ w(\cdot, T) \rightarrow \overline{u} $ in $ L^2 (0,1)$ at the same time.\\
Note that, since $ v_0 (x)\leq 0,$ 
 we deduce the following estimate
\begin{multline}\label{uT}
\|w (x,T)-\overline{u} (x)\|^2_{L^2(0,1)}=\int_0^1 \left(\int_0^T e^{v_{0} (x)
\frac{(T - \tau)}{T}} (w_{xx} (x, \tau) + f (w(x,\tau))) d \tau \right)^2 dx\\ \; \leq \;
T \parallel  w_{xx}  + f (w) \parallel^2_{L^2 (Q_T)}.
\end{multline}
{\bf Step 2:} {\it Evaluation of $  \parallel  w_{xx}  + f (w) \parallel^2_{L^2 (Q_T)}.$}  
  In this step, let us suppose that $ v_0\in C^2 ([0, 1]);$ this assumption will be removed in Step. 3.\\
Multiplying by $w_{xx}$ the equation in the first problem in \eqref{w} with $ v(x,t) = \frac{1}{T} v_{0}(x) \leq 0,$ integrating 
over $ Q_T$ and applying H\"older's inequality, we have
\begin{multline*}
\parallel  w_{xx} \parallel^2_{L^2 (Q_T)} \; = \;
\int_0^T \int_0^1 w_t w_{xx} dx dt \; - \;
\frac{1}{T} \int_0^T \int_0^1 v_{0} w w_{xx} dx dt-\int_0^T  \int_0^1   f(w)w_{xx} dx dt\\
\leq \int_0^T \int_0^1 w_t w_{xx} dx dt \; - \;
\frac{1}{T} \int_0^T \int_0^1 v_{0} w w_{xx} dx dt + \: 
\frac{1}{2}  \int_0^T  \int_0^1   f^2(w) dx dt
+ \frac{1}{2}  \int_0^T  \int_0^1   w^2_{xx} dx dt \;.
\end{multline*}
Thus, integrating by parts and recalling that $ v_0 (x) \leq 0 $, we obtain 
\begin{multline}\label{4.7}
\parallel  w_{xx} \parallel^2_{L^2 (Q_T)} \; \leq \;
2\int_0^T \int_0^1 w_t w_{xx} dx dt \; - \;
\frac{2}{T} \int_0^T  \int_0^1 v_{0} w w_{xx} dx dt 
\; + \:  \int_0^T  \int_0^1   f^2(w) dx dt\\
= \; - \int_0^T \int_0^1 (w^2_{x})_t dx dt \; + \;
\frac{2}{T} \int_0^T  \int_0^1 v_{0} w^2_{x} dx dt \; + \;
\frac{1}{T} \int_0^T \int_0^1 v_{0x} (w^2)_{x} dx dt
\; + \:   \int_0^T  \int_0^1   f^2(w) dx dt	\\
\leq \; \int_0^1  u^2_{in\,x}  dx \; - \int_0^1 w^2_{x}(T,x) dx dt  \;
 - \;
\frac{1}{T} \int_0^T  \int_0^1 v_{0xx} w^2 dx dt
\; + \:   \int_0^T  \int_0^1   f^2(w) dx dt
\\
\leq
 \int_0^1 u^2_{in\,x} 
 dx 
 + \;
\frac{1}{T} \max_{x \in [0,1]} \mid v_{0xx} \mid
\int_0^T \int_0^1   w^2 dx dt 
+ \:   \int_0^T  \int_0^1   f^2(w) dx dt.
\end{multline}
Now, we have to evaluate $ \parallel w \parallel_{C ([0,T]; L^2 (0,1))}$ and $ \parallel f(w) \parallel_{C ([0,T]; L^2 (0,1))}$. 
Since $ v_0 (x) \leq 0$, multiplying by $ w $ the equation in the first problem of (\ref{w}) and integrating by parts
 yields 
\begin{align*}
\frac{1}{2}\int_0^t \int_0^1 (w^2)_t dx ds =&\int_0^t \int_0^1 w_t w dx ds \\
=& \;
\int_0^t \int_0^1 w_{xx}w\, dx ds+
\frac{1}{T} \int_0^t \int_0^1 v_{0} w^2 dx ds+\int_0^t  \int_0^1   f(w)w dx ds
\\
\leq &
-\int_0^t \int_0^1 w^2_{x}\, dx ds+
L\int_0^T  \int_0^1   w^2 dx dt
\leq L
\int_0^T  \int_0^1   w^2 dx dt,
\end{align*}
where $ L $  is the Lipschitz constant in (\ref{1.2a}). Then, for $T\in \big(0,\frac{1}{4L}\big)$ we deduce
\begin{multline}\label{4.8}
\!\!\!\!\!\!\int_0^1 w^2 (x,t) dx \leq \int_0^1 u_{in}^2 (x) dx + 
2L\int_0^T  \int_0^1   w^2 dx dt\\
\!\leq \int_0^1 u_{in}^2 (x) dx  + 2L
T \|w\|^2_{C([0,T],L^2(0,1))}\leq  \parallel  u_{in} \parallel^2_{L^2 (0,1)}+ \frac{1}{2} \|w\|^2_{C([0,T],L^2(0,1))},\;\;t\in(0,T),
\end{multline}
so,
\begin{equation}\label{4.10}
\parallel w \parallel_{C ([0,T]; L^2 (0,1))} \; \leq \; \sqrt{2} \parallel  u_{in} \parallel_{L^2 (0,1)}.
\end{equation}
From assumption \eqref{1.2a}   and (\ref{4.10}) it follows that $f(w)\in C ([0,T]; L^2 (0,1))$ and the following estimate holds
\begin{equation}\label{4.11}
 \parallel f(w) \parallel_{C ([0,T]; L^2 (0,1))} \; \leq \; 
 L  \parallel w \parallel_{C ([0,T]; L^2 (0,1))}\leq\sqrt{2} L \parallel  u_{in} \parallel_{L^2 (0,1)}.
\end{equation}
Due to (\ref{4.7}) and (\ref{4.10})-(\ref{4.11}), we also have 
\begin{align}\label{second}
 \parallel  w_{xx}  + &f (w) \parallel^2_{L^2 (Q_T)}\leq  2\parallel  w_{xx} \parallel^2_{L^2 (Q_T)} +2\parallel  f (w) \parallel^2_{L^2 (Q_T)}\nonumber\\\; 
\leq & 2\int_0^1 u^2_{in\,x}  
dx 
 + 
\frac{2}{T} \max_{x \in [0,1]} \mid v_{0xx} \mid
\int_0^T \int_0^1   w^2 dx dt 
+ \:   3\int_0^T  \int_0^1   f^2(w) dx dt\nonumber\\
\leq& 
2\int_0^1  
u^2_{in\,x}  dx 
 + 
\frac{2}{T} \max_{x \in [0,1]} \mid v_{0xx} \mid \cdot  T  \parallel w \parallel^2_{C ([0,T]; L^2 (0,1))}
+ \:   3 T  \parallel f(w) \parallel^2_{C ([0,T]; L^2 (0,1))}
\nonumber\\
 \leq& 2\int_0^1  u^2_{in\,x} 
  dx +\big(4 \max_{x \in [0,1]} \mid v_{0xx} \mid+6TL^2 \big) \int_0^1  u^2_{in}  dx 
 \nonumber\\
 \leq &2\left(1  \; + \;
2 \max_{x \in [0,1]} \mid v_{0xx} \mid
\; + \:   3 T  L^2 \right)\, \parallel  u_{in} \parallel^2_{H^1_0 (0,1)}.
\end{align}
{\bf Step 3:} {\it Convergence of $ w (\cdot, T) $ to $\overline{u}(\cdot)$.} Note that in the previous step we can remove the assumption $ v_0\in C^2 ([0, 1]).$ Namely, if $ v_0 \not\in C^2 ([0, 1])$
we could consider  a sequence of 
uniformly bounded functions 
$\displaystyle \{ v_{0j}\}_{j \in \N},
 v_{0j} \in C^2( [0, 1]), v_{0j} (x) \leq 0, \forall x\in [0,1],$ approximating $ v_0 $ in $ L^2 (0,1).$ Then making use of the following limit relation
$$
e^{v_{0j} (x) t/T} u_{in} (x)  \mid_{t = T} \; \rightarrow  \;
e^{v_0 (x) t/T} u_{in}(x) \mid_{t = T} \; = \; \overline{u} (x) \;\;{\rm in} \;\; L^2 (0,1) \;\; {\rm as} \; j \rightarrow \infty,
$$
we conclude that equality \eqref{WT} still holds.
 Moreover,
by \eqref{uT} and \eqref{second},
we deduce that 
\begin{equation}\label{w(T)}
\|w (x,T)-\overline{u} (x)\|_{L^2(0,1)}\leq \;
C(T) \parallel  u_{in} \parallel^2_{H^1_0  (0,1)}, \textit{ where } C(T)\rightarrow0 \textit{ as } T\rightarrow0^+.
\end{equation}
{\bf Step 4:} {\it Evaluation of $ \parallel h (\cdot, T) \parallel_{L^2 (0,1)}$.}  
Multiplying by $ h$ in the equation of the second problem of (\ref{w}) and integrating by parts over $ Q_T,$ proceeding similarly to Step 2 (see in particular \eqref{4.8}) and keeping in mind (\ref{1.2a}), for every $T\in (0,\frac{1}{4L}),$ yields
\begin{align*}
\int_0^1 h^2 (x,t) dx 
\leq&\int_0^1 r_{in}^2 (x) dx
+2 \int_0^T  \int_0^1 (f(w+h) -  f(w))h \,dx dt \nonumber\\\leq &
 \int_0^1 r_{in}^2 (x) dx
+ 2L \int_0^T  \int_0^1  h^2 dx dt
\leq
 \|r_{in} \|^2_{L^2(0,1) } 
+2LT \|h\|^2_{C([0,T];L^2(0,1))}
\nonumber\\
\leq &
\|r_{in} \|^2_{L^2(0,1) } 
+\frac{1}{2}  \|h\|^2_{C([0,T];L^2(0,1))},\quad t\in (0,T).
\end{align*}
 Hence,  
\begin{equation}\label{h}
\parallel h \parallel_{C ([0,T]; L^2 (0,1))} \; \leq \; \sqrt{2} \parallel  r_{in} \parallel_{L^2 (0,1)}.
\end{equation}
{\em Conclusions.}   Thus, recalling \eqref{w}, \eqref{w(T)} and \eqref{h}, we obtain the conclusion. 
$\qquad\diamond$\\

Now we need to extend this result to the general case, 
that is, we have to prove Theorem \ref{th preserving}.

\noindent {\bf Proof (of Theorem \ref{th preserving}).} Let us fix $\eta>0$.\\
{\bf Step 0:} {\it Approximating argument. 
} Without loss of generality we can suppose $u_{in},\overline{u}\in C^1([0,1])$ and
\begin{equation}\label{u in}
|u_{in}^\prime(x_l)|=1,\qquad l=0,\ldots,n+1.
\end{equation}
Namely, if $ u_{in}, \overline{u}\not\in C^1 ([0, 1])$
we could consider two sequences in $C^1 ([0, 1]),$ 
approximating in $ L^2 (0,1),$ $ u_{in}$ and $\overline{u}, $ respectively, such that any function of the sequence approximating $u_{in}$ satisfies condition \eqref{u in}.\\
{\bf Step 1:} {\it Steering the system from $u_{in}$ to $Ku_{in},$ with $K>1.$
}  We consider the auxiliary function $\psi:[0,1]\longrightarrow\R,$ defined in the following way
$$\psi(x)=\begin{cases}
\dfrac{\overline{u} (x)}{u_{in} (x)},\quad \text{ if }\,x\in \left(0,1\right)\backslash \bigcup_{l=1}^n\left\{x_l\right\} \\
|\overline{u}^\prime (x)|, \;\;\;\text{ if }\,x=x_l,\, l=0,\ldots,n+1.
\end{cases}$$
By \eqref{u in} we have that $\displaystyle\lim_{x\rightarrow x_l}\frac{\overline{u} (x)}{u_{in} (x)}=|\overline{u}^\prime (x_l)|,\;\; \forall l=0,\ldots,n+1,$
so $\psi\in C([0,1]).$ Let us introduce the universal constant $$\displaystyle K=K(u_{in},\overline{u}):=
\max_{x\in [0,1]}\psi(x)+1>1
.$$
For any $\displaystyle 0 < \rho <\frac{\rho_0}{2}:= \frac{1}{2}\min_{l=0,\ldots,n} \big\{ x_{l+1}-x_{l}\big\},$ consider
the following set \\
$$\displaystyle A_{\rho}:=
\bigcup_{l=0}^{n} \Big(x_l + \rho, x_{l+1} -\rho\Big).
$$
From the definition of $K$ it follows that
\begin{equation}\label{5.1}
\displaystyle K
>\max_{x\in\overline{A}_{\rho}}
\Big\{\frac{\overline{u} (x)}{u_{in} (x)}\Big\},
\;\;\;
\;\forall \rho\in \left(0,\frac{\rho_0}{2}\right).
\end{equation} 
Let us select 
\begin{equation*}
v (x,t) = m := \frac{\ln K}{t_1}, \;\; (t,x) \in (0,1)\times(0, t_1),
\end{equation*} 
for some arbitrarily small $ t_1>0$. Then, let us apply the auxiliary constant bilinear control  $ v (x,t) = m > 0, \forall x\in (0,1)$ on  the interval $(0, t_1).$ 
For $t=t_1$ the solution of the first problem in \eqref{w} has the following representation in Fourier series
\begin{align}\label{w conv}
w (x,t_1) &= \;  H (x, t_1) \; + \; e^{mt_1} \sum_{p = 1}^\infty  2 e^{-(p \pi )^2 t_1} \left(\int_0^1 u_{in}(r) \sin p\pi r \, dr \right) \sin p\pi  x\nonumber\\
&= \; H (x, t_1) \; + \; e^{mt_1} \sum_{p = 1}^\infty  2 (e^{-(p \pi)^2 t_1} -1) \left(\int_0^1 u_{in} (r) \sin p\pi r \, dr \right) \sin p\pi  x
\; + \; e^{mt_1} u_{in} (x)\nonumber\\
&= \; H (x, t_1) \; + R(x,t_1)
+ K u_{in} (x),
\end{align}
where 
\begin{align*}
H (x,t_1) & := \; \sum_{p = 1}^\infty  2 \left[\int_0^{t_1} e^{(m -(p \pi )^2 ) (t_1 -t)} \left(\int_0^1 f (w (r, t))  \sin p\pi r \, dr \right) dt \right]\sin p\pi  x,\nonumber\\
R (x,t_1) & :=
K \sum_{p = 1}^\infty  2 (e^{-(p \pi)^2 t_1} -1) \left(\int_0^1 u_{in} (r) \sin p\pi r \, dr \right) \sin p\pi  x\,.
\end{align*}
By Parseval's equality, we deduce 
\begin{align}\label{R}
\!\!\!\!\!\|R (\cdot,t_1)\|^2_{L^2(0,1)}\!\!\!&= K^2\left\| \sum_{p = 1}^\infty  (e^{-p\pi^2 t_1} -1) \left(\int_0^1 u_{in} (r) \sqrt{2}\sin p\pi r \, dr \right) \sqrt{2}\sin p\pi  x\,\right\|^2_{L^2(0,1)}\nonumber\\
&= K^2
\sum_{p = 1}^\infty(e^{-p\pi^2 t_1}-1)^2 \left|\int_0^1 u_{in} (r) \sqrt{2}\sin p\pi r \, dr  
\right|^2
\nonumber\\
&\leq K^2(1-e^{-\pi^2 t_1})^2
\sum_{p = 1}^\infty \left|\int_0^1 u_{in} (r) \sqrt{2}\sin p\pi r \, dr  
\right|^2
\nonumber\\
&= K^2(1-e^{-\pi^2 t_1})^2\|u_{in}\|^2_{L^2(0,1)}.
\end{align}
In the same
way, 
using assumption \eqref{1.2a}, $f(0)=0,$
 and H\"older's inequality we obtain
\begin{align}\label{H}
\parallel H (\cdot, t_1) \parallel^2_{L^2 (0,1)} &=
\left\|\sum_{p = 1}^\infty   
\left[\int_0^{t_1} e^{(m -(p \pi )^2 ) (t_1 -t)}\left( \int_0^1  
f (w (r, t))  \sqrt{2}\sin p\pi r \, dr \right)
 dt \right]\sqrt{2}\sin p\pi  x\right\|^2_{L^2 (0,1)}
\nonumber\\ 
&=  \,
 \sum_{p = 1}^\infty \left| \int_0^{t_1} e^{(m -(p \pi )^2 ) (t_1 -t)}\left( \int_0^1  
f (w (r, t))  \sqrt{2}\sin p\pi r \, dr \right)
 dt
\right|^2
\nonumber\\
&\leq  \,
 \sum_{p = 1}^\infty \left( \int_0^{t_1}  e^{2(m -(p \pi )^2 ) (t_1 -t)}dt
\right)\int_0^{t_1}\left|\int_0^1f (w (r, t))  \sqrt{2}\sin p\pi r \, dr\right|^2\!dt 
\nonumber\\
 &\leq
 \sum_{p = 1}^\infty 
 e^{2mt_1}t_1
\int_0^{t_1}\left|\int_0^1f (w (r, t))  \sqrt{2}\sin p\pi r \, dr\right|^2\!dt 
\nonumber\\
&=  \,
K^2 
t_1 
\int_0^{t_1}\sum_{p = 1}^\infty\left|\int_0^1f (w (r, t))  \sqrt{2}\sin p\pi r \, dr\right|^2\!dt 
=K^2 
t_1 
\int_0^{t_1}
\int_0^1f ^2(w (r, t)) drdt 
\nonumber\\
&
\leq L^2\,K^2 
t_1
\int_0^{t_1}\int_0^1
 w^2 (r, t) dr
\,dt
\leq L^2\,K^2 t^2_1 \parallel  w \parallel^2_{C([0,t_1]; L^2 (0,1))}.
\end{align}
Now, we have to evaluate $ \| w\|_{C ([0,t_1]; L^2 (0,1))}.$  
Multiplying by $ w $ the equation in the first problem of (\ref{w}), integrating by parts, and arguing as in the proof of
\eqref{4.8}, for every $t\in (0,t_1),$
we have
\begin{align*}
\frac{1}{2}\int_0^t \int_0^1 (w^2)_t dx ds 
&=
-\int_0^t \int_0^1 w^2_{x}\, dx ds+m \int_0^t \int_0^1 w^2 dx ds+\int_0^t  \int_0^1   f(w)w dx ds\\
&\leq (m+L)
\int_0^t  \int_0^1   w^2 dx ds,
\end{align*}
where $ L $  is as in (\ref{1.2a}).
Then
$$
\int_0^1 w^2(x,t) dx 
\leq \int_0^1 u_{in}^2(x) dx 
+ 2(m+L)
\int_0^t  \int_0^1   w^2 dx ds,\;\;t\in (0,t_1).
$$
Thus, applying Gr\"onwall's inequality we deduce 
$
\|w(t,\cdot)\|^2_{L^2(0,1)}
\leq e^{2(m+L)t_1}\|u_{in}\|^2_{L^2(0,1)},\;t\in (0,t_1).
$
So,
\begin{equation}\label{w energy}
\|w\|_{C([0,t_1];L^2(0,1))}
\leq e^{(m+L)t_1}\|u_{in}\|_{L^2(0,1)}=Ke^{L\,t_1}\|u_{in}\|_{L^2(0,1)}\;\; .
\end{equation}
Making use of \eqref{w conv}-\eqref{w energy}, we have that
\begin{align}\label{w final} \|w(\cdot,t_1)-K u_{in}\|_{L^2(0,1)}&=
\|H (x, t_1) \; + R(x,t_1)\|_{L^2(0,1)}\nonumber\\
&\leq K\left[(1-e^{-\pi^2 t_1})+t_1K\,L\,e^{L\,t_1}\right]\|u_{in}\|_{L^2(0,1)}.
\end{align}
Now, we evaluate $ \parallel h (\cdot, t_1) \parallel_{L^2 (0,1)}$. 
Multiplying by $ h$ both members of
 the equation in the second problem of (\ref{w}) and integrating by parts, proceeding similarly to \eqref{w energy} and to the proof of Lemma \ref{lem preserving}, and keeping in mind (\ref{1.2a}) we obtain
\begin{align*}
\int_0^1 h^2 (x,t) dx 
&\leq\int_0^1 r_{in}^2 (x) dx
+2\,m \int_0^{t_1}  \int_0^1  h^2 dx dt+2 \int_0^{t_1}  \int_0^1 (f(w+h) -  f(w))h \,dx dt \nonumber\\&\leq 
 \int_0^1 r_{in}^2 (x) dx
+ 2(m+L) \int_0^{t_1}  \int_0^1  h^2 dx dt\,.
\end{align*}
 Hence, using Gr\"onwall's inequality, we have 
 $\|h(t,\cdot)\|^2_{L^2(0,1)}
\leq e^{2(m+L){t_1}}\|r_{in}\|^2_{L^2(0,1)},\;t\in (0,t_1),$ so
 \begin{equation}\label{h energy}
\|h(t_1,\cdot)\|^2_{L^2(0,1)}\leq\|h\|_{C([0,t_1];L^2(0,1))}
\leq 
Ke^{L\,t_1}\|r_{in}\|_{L^2(0,1)}.
\end{equation}
\noindent 
Thus, by \eqref{w final} and \eqref{h energy}, there exists $t_1=t_1(\eta)>0,\,t_1\ll1$ such that the following inequality holds
\begin{equation}\label{5.10}
\|u (\cdot, t_1) - K u_{in} (\cdot)\|_{L^2(0,1)}\leq \frac{\sqrt{2}}{8}\eta +  Ke^{L}\|r_{in}\|_{L^2(0,1)}.
\end{equation}
{\bf Step 2:} {\it Steering the system from $K u_{in}+r_{in}
$ to $\overline{u}.$} 
In this step, let us represent again the solution $ u $ of \eqref{u pres proof} as a sum of two functions $ w (x,t) $ and $h (x,t)$, which solve the problems in \eqref{w} in $(0,1)\times(t_1,T),$ with
 the modified initial states $K u_{in} $ instead of $ u_{in} $ and $r_{in}(\cdot)=u (\cdot, t_1) - K u_{in} (\cdot).$
\\
By \eqref{5.1} it follows that
\begin{equation}\label{ass step}
\frac{\overline{u} (x)}{K u_{in} (x)}<1, \;\;\qquad \forall x \in A_{\rho},\, 
\forall \rho\in \left(0,\frac{\rho_0}{2}\right),
\end{equation} 
moreover,
\begin{equation}\label{nu}
\forall \rho \in \left(0,\frac{\rho_0}{2}\right),\; \exists\, \nu=\nu(\rho)>0 \;:
\; \;\nu\leq\frac{\overline{u} (x)}{K u_{in} (x)}, \;\; \forall x \in A_{\rho}.
\end{equation} 
Keeping in mind the proof of Lemma \ref{lem preserving}, we consider
the following function defined on $[0,1]$
\begin{equation}\label{5.11}
v_0 (x) \; = \;  \left\{ \begin{array}{ll}
\ln \left( \frac{\overline{u} (x)}{K u_{in} (x) } \right), \;\;\;\; &  x \in A_{\rho},\\
0, \;\;\;\; &  {\rm elsewhere} \; {\rm in} \; [0, 1]\,. \\
\end{array}
\right.   
\end{equation}
By \eqref{nu} 
$ v_0 \in L^\infty (0, 1),$ and
by \eqref{ass step} we have 
$v_0 (x) \; \leq 0,\;\forall x\in [0,1].$
Let us select the bilinear control
$$
v (x,t) \; := \; \frac{1}{T} v_0 (x)\;\;\;\forall x\in(0,1)\times(t_1,T).
$$
We remark that, thanks to \eqref{ass step}-\eqref{nu}, the assumption 
\eqref{ass} of Lemma \ref{lem preserving} holds. Then proceeding similarly to Lemma \ref{lem preserving}, with a proof essentially identical, there  exists $T=T(\eta)>t_1$ with $T-t_1$ sufficiently small ($0<T-t_1\ll 
\frac{1}{4L}
$)
 we can obtain the following inequality (similar to \eqref{ineq pres} of Lemma \ref{lem preserving}) 
\begin{equation}\label{5.12}
\|u (\cdot, T) - \overline{u}_{\rho} \|_{L^2(0,1)}\leq \frac{\eta}{4}+\sqrt{2}\|u (\cdot, t_1) - K u_{in} (\cdot)\|_{L^2(0,1)},
\end{equation}
with
$$
\overline{u}_\rho(x)\; = \;  \left\{ \begin{array}{ll}
\overline{u} (x), \;\;\;\; &  x \in A_\rho,\\ 
0, \;\;\;\; &  {\rm elsewhere} \; {\rm in} \; [0, 1].\\
\end{array}
\right.   
$$  
Note that here exists $
\overline{\rho}=\overline{\rho}(\eta)>0$ such that
\begin{equation}\label{approx u}
\displaystyle 
\parallel  \overline{u}_{\overline{\rho}} -\overline{u} \parallel_{L^2 (0,1)}<\frac{\eta}{2}.
\end{equation}
Then, from \eqref{5.12}, \eqref{approx u} and \eqref{5.10} we obtain the conclusion
$$
\parallel  u(\cdot,T) -\overline{u}(\cdot) \parallel_{L^2 (0,1)}\leq \parallel  u(\cdot,T) - \overline{u}_{\overline{\rho}}(\cdot) \parallel_{L^2 (0,1)}+\parallel   \overline{u}_{\overline{\rho}} -\overline{u} \parallel_{L^2 (0,1)}
\leq\eta+ \sqrt{2}Ke^{L}\|r_{in}\|_{L^2(0,1)}\,.\quad\diamond
$$

\end{document}